\documentclass[12pt]{amsart}
\usepackage{RWF}
\usepackage{pagesize}

\begin{document}

\title{On the Rozansky-Witten weight systems}
\author{Justin Roberts}
\address{Department of Mathematics, UC San Diego, 9500 Gilman Drive,
La Jolla CA 92093, USA}
\email{justin@math.ucsd.edu}
\author{Simon Willerton}
\address{Department of Pure Mathematics, University of Sheffield,
Hicks Building, Sheffield, S3 7RH, UK}
\email{s.willerton@sheffield.ac.uk}
\date{Feb 27, 2006}


\begin{abstract}
Ideas of Rozansky and Witten, as developed by Kapranov, show that
a complex symplectic manifold $X$ gives rise to Vassiliev weight
systems. In this paper we study these weight systems by using
$D(X)$, the derived category of coherent sheaves on $X$. The main
idea (stated here a little imprecisely) is that $D(X)$ is the
category of modules over the shifted tangent sheaf, which is a Lie
algebra object in $D(X)$; the weight systems then arise from this
Lie algebra in a standard way. The other main results are a
description of the symmetric algebra, universal enveloping
algebra, and Duflo isomorphism in this context, and the fact that
a slight modification of $D(X)$ has the structure of a braided
ribbon category, which gives another way to look at the associated
invariants of links. Our original motivation for this work was to
try to gain insight into the Jacobi diagram algebras used in
Vassiliev theory by looking at them in a new light, but there are
other potential applications, in particular to the rigorous
construction of the $(1+1+1)$-dimensional Rozansky-Witten TQFT,
and to hyperk\"ahler geometry.
\end{abstract}

\maketitle
\tableofcontents

\pg


\section*{Introduction}

\subsection*{Motivation.}

The {\em Kontsevich integral} is a beautiful and very powerful
invariant of framed knots in $S^3$. It takes values in a certain
graded algebra $\A$ of {\em Jacobi diagrams}, and is universal for the
class of {\em Vassiliev {\rm (}finite-type{\rm )} invariants}, as well
as determining all the {\em quantum invariants} (Jones polynomial
et. al.)  associated to quantum groups. The definitive exposition is
by Bar-Natan \cite{BN}.

Over the last few years the theory of the Kontsevich integral has been
considerably extended (see Le, Murakami and Ohtsuki \cite{LMO} and
Murakami and Ohtsuki \cite{MurTQFT}), resulting in a system of
Kontsevich-like invariants for links, graphs, and $3$-manifolds
possessing much of the functoriality of a topological quantum field
theory. Despite these successes, basic questions about the topological
interpretations of the Kontsevich invariant and of the algebra $\A$
itself remain largely unanswered.

The standard way to study diagram spaces such as $\A$ is by means
of {\em weight systems}, ie.\ functions on it, which are most
easily obtained from Lie algebras. A finite-dimensional Lie
algebra $\g$ with an invariant non-degenerate metric defines
weight system homomorphisms from diagram spaces to spaces of {\em
invariant tensors} on $\g$; recognisable formulae often emerge
from looking in this way at diagrammatic identities ``at the level
of Lie algebras'', as for example in Bar-Natan, Garoufalidis,
Rozansky and Thurston \cite{BGRT}. From this has emerged the idea,
pursued in particular by Vogel \cite{Vogel}, that the diagrams
themselves form some kind of universal Lie-algebra-like structure.

In this paper we propose to study diagram algebras from an alternative
point of view using {\em Rozansky-Witten weight systems}
\cite{RW}. These arise from {\em complex symplectic manifolds},
according to Kapranov and Kontsevich \cite{Kapranov, KontsRW}, and map
diagram algebras to {\em Dolbeault cohomology groups} of such
manifolds. Our original motivation for this study was to try to
understand the extent to which diagrams behave like elements of
cohomology; we were seeking to interpret $\A$ as some kind of ring of
universal characteristic classes, and had already been studying
certain diagrammatic formulae as if they were cohomological
identities.

This is a reasonable point of view: after all, Kontsevich's
formulation of {\em graph cohomology} \cite{KontsevichFeynman} shows
that indeed, diagrams may be thought of as representing elements of
homology and cohomology, though this combinatorial framework affords
little topological insight. Although Kontsevich has also given
interpretations of graph cohomology via Gelfand-Fuchs cohomology and
Lie algebras of formal Hamiltonian vector fields, we still hope that
there is a more direct explanation for much of the theory. We would
like to be able to view graph cohomology as the cohomology of some
kind of interesting and meaningful geometrical classifying space (by
analogy with {\em fatgraph} cohomology, which is the cohomology of the
moduli space of Riemann surfaces), and then use the geometry of this
space to give new explanations of the existence and properties of the
knot and $3$-manifold invariants. There is an obvious candidate, {\em
outer space} \cite{KontsevichFeynman}, but it still seems rather too
abstract for our purposes, and these goals remain
unfulfilled. Fortunately, Rozansky-Witten theory is a fruitful subject
to study in its own right.

In this paper we deal only with the nature of the Rozansky-Witten
weight systems. That is, we are looking at diagrams ``at the level of
complex symplectic manifolds'', and studying the analogies between Lie
algebra and Rozansky-Witten weight systems, in a sense parallelling
the paper \cite{BGRT}. An alternative focus would be to use the theory
to derive results about hyperk\"ahler geometry, in the manner of
Hitchin and Sawon \cite{Hitchin-Sawon}, but we will avoid this here.
Likewise, though we touch here on the Rozansky-Witten link invariants,
we will for the most part postpone the study of the associated
topological invariants for a sequel in which we apply our techniques
to set up the full $(1+1+1)$-dimensional Rozansky-Witten TQFT.

It is on the face of it very surprising that objects as disparate
as Lie algebras and complex symplectic manifolds give rise to
weight systems. The main point of our paper is to unify these two
worlds, showing how to define and handle Rozansky-Witten weight
systems in a way completely analogous to the Lie algebra ones. We
show in fact that a complex symplectic manifold gives rise to
something akin to a metric Lie algebra, and then investigate the
ramifications of this analogy. The catch here is that this
something is an object in a category other than the usual category
of vector spaces; in fact, the category must be taken to be the
{\em bounded derived category of coherent sheaves} on the
manifold.

Now phrases like this used to strike terror into the hearts of the
authors, and we presume some readers will also recoil slightly!  But
we are really convinced that the use of derived categories gives the
most natural and elegant formulation of the Rozansky-Witten theory,
and hope by our exposition to convince the reader likewise. An
additional justification for their use is that in our forthcoming
construction of the $(1+1+1)$-dimensional Rozansky-Witten TQFT, the
derived category is essentially the category associated to a circle.

It is also worth mentioning here a disadvantage of our approach, which
is that the beautiful $L_\infty$ structure described by Kapranov is
thrown away. We would need this if we were interested in weight
systems defined on higher graph cohomology, but these do not figure in
the computation of the usual knot or TQFT invariants, so this is an
acceptable loss.

A brief sketch of the results in this paper appears in \cite{R}, and
further questions appear in \cite{RS}.


\subsection*{Overview of results}
\subsubsection*{The derived category as a Lie algebra representation category}

Kapranov showed that one could use a certain $L_\infty$--algebra
structure to obtain the Rozansky-Witten weight systems: our approach
is to work in the derived category and use Lie algebras type objects
as in Chern-Simons theory.  The first step is to identify the derived
category as the representation category of a certain Lie algebra
inside it.
\begin{introthm}
Let $X$ be a complex manifold.  The shifted tangent bundle $T[-1]$
is a Lie algebra object in the bounded derived category $D(X)$,
and $D(X)$ is the category of modules over $T[-1]$.
\end{introthm}

To explain this, first note that a Lie algebra object in a additive
symmetric tensor category means an object $L$ in the category with a
bracket morphism $L\otimes L\to L$ which satisfies suitable versions
of the Jacobi and anti-symmetry relations, so that a Lie algebra
object in the category of vector spaces is a usual Lie algebra and a
Lie algebra object in the category of graded vector spaces is a graded
Lie algebra.  A module $M$ over $L$ is then an object of the category
with an action morphism $M\otimes L\to M$ satisfying an appropriate
condition.  Next we need to know a little about the derived category.
This has as objects bounded chains complexes of coherent sheaves on
$X$, so in particular for each coherent sheaf $\scrE$ on $X$ there is
the $\scrE[-i]$, consisting of $\scrE$ in degree $i$ and zero
elsewhere (we still write $\scrE$ for $\scrE[0]$).  One fundamental
fact about the derived category is that hom-sets can be identified
with $\Ext$ groups, or equivalently, cohomology groups, so that
\[\Mor_{D(X)}(\scrE,\scrF[i])\cong\Ext^i(\scrE,\scrF).\]
The next key ingredient is the Atiyah class $\alpha_\scrE$ for a
coherent sheaf $\scrE$; this is a characteristic class which lives
in $\Ext^1(\scrE, \scrE \otimes \Omega)$, which we can identify as the
hom-set $\Hom_{D(X)}(\scrE\otimes T[-1],\scrE)$. 

Thus the Atiyah class can be thought of as a morphism
$\alpha_\scrE\colon \scrE\otimes T[-1]\to \scrE$. In particular the
Atiyah class of the tangent bundle gives the Lie bracket $T[-1]\otimes
T[-1]\to T[-1]$, and the other Atiyah classes $\alpha_\scrE$ give
module maps.

Unfortunately these do not give the action of $T[-1]$ on every object
of the derived category. So in fact we take the more elegant approach
of realising the action of $T[-1]$ as a natural transformation
$\alpha$ from $\id\otimes T[-1]$ to the identity functor of
$D(X)$. This gives for every object $A$ in the derived category a
morphism $\alpha_A\colon A\otimes T[-1]\to A$, and naturality ensures
that every morphism $A\to B$ intertwines the action on $A$ and $B$.

The natural transformation $\alpha\colon \id \otimes T[-1]\to\id$ is
obtained using an {\em integral transform}. It is a standard principle
of ``correspondences'' that objects of $D(X\times X)$ define functors
$D(X) \ra D(X)$, and that morphisms of $D(X \times X)$ define natural
transformations between them: in fact we get a functor from $D(X\times
X)$ to the functor category $\Fun(D(X),D(X))$. Our natural
transformation $\alpha$ is obtained from a morphism $\Odel \otimes
\pi^*T[-1]\to \Odel$ in $D(X\times X)$ which is essentially one half
of the Atiyah class of the structure sheaf of the diagonal.

\subsubsection*{Complex symplectic manifolds and invariant metrics}
A complex symplectic form is a non-degenerate holomorphic two-form, in
other words an element of $H^0(X, \bigwedge^2T^\dualstar)$ and we can
identify this with a symmetric element of the hom-set
$\Mor_{D(X)}(T[-1]\otimes T[-1], \OO_X[-2])$.  This isn't quite an
invariant metric on the Lie algebra object $L=T[-1]$: such a thing
would be a symmetric $L$-module map $L\otimes L\to \one$, but we have
an extra shift $[-2]$.  To handle this we can work in the ``extended
derived category'' $\tilde D(X)$ whose hom-set $\Hom_{\tilde D(X)}(A,
B)$ is the graded group $\Ext^\gradestar(A,B)$, where the shift
problem disappears. Thus 
\begin{introthm}
Let $X$ be a complex symplectic manifold.  The shifted tangent bundle
$T[-1]$ is a metric Lie algebra object in the extended bounded derived
category $\tilde D(X)$, and $D(X)$ is the category of modules over
$T[-1]$.
\end{introthm}

\subsubsection*{Symmetric and universal enveloping algebras}
The reader familiar with ideas of Vassiliev invariants will know that
other Lie algebraic concepts such as the universal enveloping algebra,
symmetric algebra, Poincar\'e-Birkoff-Witt isomorphism and Duflo
isomorphism play important roles. We show that analogues of these
makes sense for the Lie algebra object $T[-1]$ on any complex
manifold. If we were working in an abelian category then the symmetric
and universal enveloping algebras could be constructed as quotients or
subobjects of the tensor algebra of the Lie algebra object, but as the
derived category is not abelian, we have to work a bit harder.  The
symmetric algebra of $T[-1]$ is easily identifiable as
$S=\bigoplus(\bigwedge^iT)[-i]$, the shifted exterior algebra on
$T$. But the universal enveloping algebra is less obvious: we define
$U=\pi_*\sHom(\Odel,\Odel)$, and prove that $D(X)$ can be thought of
as the representation category of $U$ also:

\begin{introthm}
The object $U\in D(X)$ is an associative algebra object. There is a
canonical map $L=T[-1]\to U$ with respect to which $U$ is the
universal enveloping algebra of $L$. The algebra $U$ acts on all
objects of $D(X)$ in a manner compatible with the action of $T[-1]$.
\end{introthm}

The Poincar\'e-Birkhoff-Witt and Duflo isomorphisms have their
analogues in this world. For standard Lie algebras, the natural
symmetrisation map $\PBW: S(\g) \ra U(\g)$ is an isomorphism of
$\g$-modules and hence induces an isomorphism on their invariant
parts. The latter can be corrected by a strange automorphism of
$S(\g)^\g$ to give Duflo's algebra isomorphism between $S(\g)^\g$ and
$Z(\g) =U(\g)^\g$. 

In our context, there is again an isomorphism of objects $\PBW: S
\cong U$ in $D(X)$. Our proof is an elaboration of ideas of Markarian
\cite{Ma}. The correct categorical way to think of invariants is as
homomorphisms from the trivial object, which amounts in the derived
category to cohomology. Thus the induced map on invariant parts gives
an isomorphism between the polyvector field cohomology
\[ \HT^*(X) = H^*(S) = H^\gradestar(X,\textstyle\bigwedge^\gradestar T)\]
and the Hochschild cohomology
\[ \HH^*(X) =  H^*(U) = \Ext^*(\Odel, \Odel)\]
first demonstrated by Gerstenhaber and Schack. The analogue of Duflo's
isomorphism between these algebras is Kontsevich's ``theorem on a
complex manifold''. Although this isomorphism exists for all complex
manifolds, in the case of complex symplectic manifolds it follows from
the wheeling theorem of Bar-Natan, Le and Thurston \cite{BLT}.

A corollary of these theorems is the existence, given a complex
symplectic manifold $X$, of sheaf-cohomology-valued Vassiliev weight
systems defined on all the usual algebras of Jacobi diagrams,
naturally compatible with operations such as gluing of legs, etc.

\subsubsection*{Ribbon categories and link invariants}

The theory of the Knizhnik-Zamolodchikov equation gives a way to
produce an interesting ribbon category structure on the category of
representations of $U\g \otimes \C[[\hbar]]$, which by Drinfeld's work is
equivalent to the category of representations of a quantum group. 
This result has an analogue in our context. 

\begin{introthm} 
The category $\tilde D(X)$ has a natural non-symmetric ribbon tensor
category structure when $X$ is a complex symplectic manifold. 
\end{introthm}

Ribbon categories automatically define framed link invariants. The
ones arising from $\tilde D(X)$ agree with the invariants obtained by
taking the Kontsevich integral and composing with the weight systems;
they may be thought of as the ``knot polynomial'' type quantum
invariants arising from complex symplectic manifolds. We do not
however know of any analogue of Drinfeld's theorem in this context.


\subsection*{Outline of the paper}

The first two sections of the paper are an exposition of the
``standard'' approach to Rozansky-Witten weight systems.

In Section 1 we give a brief description of what weight systems are,
and how they are obtained from finite-dimensional metric Lie algebras.

In Section 2 we describe, by analogy with Chern-Weil theory, the
differential-geometric formulation of Rozansky-Witten invariants as
integrals of suitable curvature forms. The original treatment of
Rozansky and Witten used physics (path integrals) as a motivation and
Riemannian geometry for the actual construction of weight systems for
hyperk\"ahler manifolds. We follow instead Kapranov's reworking in
terms of hermitian differential geometry, which has the advantage of
demonstrating that the construction does not actually depend on the
hyperk\"ahler metric, and will work for any complex symplectic
manifold.

The next sections are essentially reformulations of the first two,
introducing the language in which our theorems are going to be stated.

In Section 3 we reformulate the construction of weight systems from
metric Lie algebras so that it generalises to metric Lie algebras in
categories other than the category of vector spaces. This is all based
on work of Vogel and Vaintrob \cite{Vogel,Vaintrob}.

In Section 4 we explain the language of derived categories (first in a
general way and then with specific reference to sheaf theory), which
will be necessary in Section 5 when we reformulate the relevant
differential geometry in terms of sheaf theory, following
Kapranov. The key concept is the {\em Atiyah class}, the cohomological
version of the curvature of a holomorphic bundle. In Section 6 
we show how it gives a Lie bracket.

In Section 7 we explain various generalisations of the idea of a
weight system to other graph algebras, and how these relate to
Lie-theoretic concepts such as symmetric and universal enveloping
algebras. In Section 8 we show how these concepts manifest themselves
in the context of complex symplectic manifolds and how they give more
interesting kinds of weight systems. 

In Section 9 we show how to turn $D(X)$ into a ribbon category,
thereby giving another way to explain the associated invariants of
links. The paper concludes in Section 10 with a summary of the analogy
between the world of Lie algebras and complex symplectic manifolds,
which we will extend in the sequel to an analogy between Chern-Simons
TQFT and Rozansky-Witten TQFT.


\begin{ack}
The first author was partially supported by EPSRC, the NSF, and a
JSPS research fellowship at RIMS, Kyoto. The second author was
partially supported by EPSRC, the NSF, UCSD project for geometry
and physics, a Marie Curie fellowship from the European Union and
the Department of Social Security. We are indebted to Tom
Bridgeland for teaching us about derived categories, and to Justin
Sawon, Alexei Bondal, Mikhail Kapranov, Thang Le, Nikita
Markarian, Boris Shoikhet and Arkady Vaintrob for various helpful
discussions. We apologise for the very long delay in finishing this
paper. 
\end{ack}

\pg


\section{Lie algebra weight systems}
\label{Section:LAweightsystems}
 We begin with a brief description
of the algebra $\A$ of Jacobi diagrams used in Vassiliev theory,
and of how it is studied using weight systems arising from
finite-dimensional metric Lie algebras. The more involved parts of
the theory are deferred until
Section~\ref{Section:LAweightsystemsAgain}. Apart from our grading
convention, this is all standard, see Bar-Natan \cite{BN}.


\subsection{Jacobi diagrams}

The Kontsevich integral is an invariant of framed oriented knots
in $S^3$. It takes values in the complex, graded, algebra $\A$ of
{\em Jacobi diagrams} defined as follows. Consider all isomorphism
classes of connected trivalent graphs containing a preferred
oriented circle and with a choice of cyclic orientation at each
vertex not on the preferred circle. (The ones on the circle are
canonically oriented because the circle is oriented.) Define $\A$
to be the complex span of such classes, quotiented by the
vertex-antisymmetry, and IHX relations, pictured below.  When the
IHX relation involves an edge in the preferred circle it is called
the STU relation.  In this paper we will grade $\A$ by the total
number of vertices of the graph, which is even.  It is important
to note that this is {\em twice} the conventional grading.

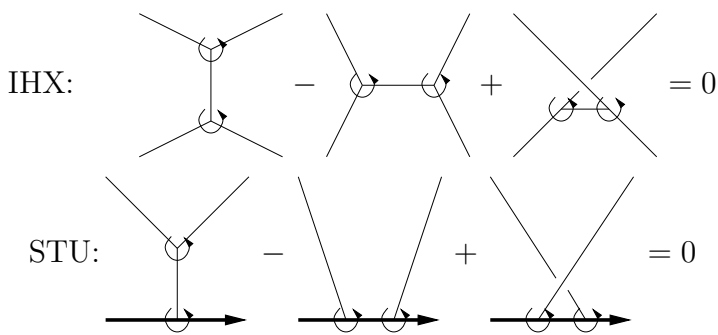
\[ \text{IHX:}\qquad\vpic{I} - \vpic{H} + \vpic{X} =0\]
\[ \text{STU:}\vpic{S}-\vpic{T} +\vpic{U} =0\]

Such graphs are usually described using planar pictures, in which the
preferred circle is drawn as an external loop, and the rest of the
graph is inscribed. One thinks of it as a graph ``with legs'' which is
attached to the external circle. One useful function of such a planar
projection is that every vertex may be given the canonical
``anticlockwise'' orientation, so the orientations need not be drawn
explicitly. Any oriented abstract graph may be drawn in such a
way. The antisymmetry relation then can be denoted by the following
picture:

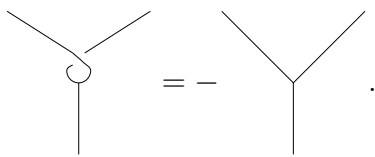
\[ \vpic{vxL} = - \vpic{vxR}.\]

The space $\A$ is a commutative algebra, whose product $\#$ is given
by connect-summing diagrams arbitrarily along their preferred oriented
circles.


\subsection{Lie algebra weight systems}
\label{Subsection:LAweightsystems}
 To obtain numerical knot
invariants from the Kontsevich integral, or simply to study the
infinite-dimensional space $\A$, it is necessary to construct
linear maps from $\A$ to some better understood rational vector
space such as $\Q$. (In fact we will typically work with complex
vector spaces in this paper.) Such a map is called a {\em weight
system}.

The simplest way to obtain a weight system taking values in $\C$ is to
pick a finite-dimensional Lie algebra $\g$ with a metric $b$ (a
non-degenerate invariant symmetric bilinear form), and a
finite-dimensional representation $V$ of $\g$. This information is
completely encoded by the following three $\g$-module maps:
\[   a=[-,-]:  \g \otimes \g \ra \g \qquad b: \g \otimes \g \ra \C
\qquad a_V: V \otimes \g \ra V. \]
Since $b$ induces an isomorphism $\g \cong \g^\dualstar$, we may
rewrite the Lie bracket as a skew trilinear form and think of it as a
tensor $f \in \bigwedge^3\g^\dualstar$.  Additionally, we may ``invert''
the metric to define a Casimir element $c \in S^2\g$. The action $a_V$
is usually thought of as a tensor in $V^\dualstar \otimes V \otimes
\g^\dualstar$.

Now, a graph in $\A$ defines a way of contracting together these
tensors to obtain a scalar in $\C$. Simply insert $f$ at the internal
vertices, $a_V$ at the external vertices, and $c$ on the internal
edges, and contract the $\g$-$\g^\dualstar$ pairs and the
$V$-$V^\dualstar$ pairs as indicated by the graph. The
vertex-orientation corresponds precisely to the information needed to
insert $f$ at a vertex; without it there would be a sign
ambiguity. The symmetry of $c$ means that no orientation on the edge
is required. It is easy now to check that the relations in $\A$ are
satisfied by this assignment, and that the weight system $w_{\g, V}:
\A \ra \C$ is well-defined.


\pg

\section{Rozansky-Witten weight systems}
\label{Section:RWweightsystems}
 In this section we explain,
following Kapranov \cite{Kapranov}, a construction via hermitian
differential geometry of weight systems from complex symplectic
manifolds. We are not actually going to use this approach in the
rest of the paper, but it's likely that at first sight it will be
more illuminating than the later sheaf cohomology approach; in any
case, Kapranov's paper is a little terse, and we feel it is
worthwhile to expand on his construction. Actually, his
demonstration of the Lie structure only works for K\"ahler
manifolds, so by extending this to all complex manifolds we are
tidying up a little too.

\subsection{Chern-Weil theory}

In this context, it is natural to consider Rozansky-Witten theory as
an variant of Chern-Weil theory. Instead of using the curvature of a
smooth connection on a smooth complex vector bundle to give invariants
in the de Rham cohomology of the base manifold, we will use a the
curvature of a hermitian connection on a holomorphic vector bundle to
give invariants in the Dolbeault cohomology of the base complex
symplectic manifold.

Recall that if E is a smooth complex vector bundle over the smooth
manifold $X$ and $\Omega^p(X;E)$ is the space of smooth $p$-forms with
values in $E$, then there is no canonical choice of differential on
$\Omega^\gradestar(X;E)$. But if we pick a smooth connection on $E$,
that is a covariant derivative $\nabla: \Omega^0(X;E) \ra
\Omega^1(X;E)$, then we induce operators $\nabla: \Omega^p(X;E) \ra
\Omega^{p+1}(X;E)$.  It is a standard fact that the composite
$\nabla^2$ is given by wedging with the {\em curvature two-form} $F
\in \Omega^2(X; \End(E))$. One can then use $GL(E)$-invariant
polynomials in $F$ to define cohomology classes which are independent
of the choice of connection.  When the bundle $E$ has rank $r$ these
polynomials are spanned by the functions $F \mapsto \tr (F^d)$, for $0
\leq d\leq r$.  The resulting cohomological invariants of $E$,
suitably normalised, define its Chern classes modulo torsion; more
precisely, the class
\[\ch_d(E) = \left[\frac{1}{d!}\tr \left(\frac{-F}{2 \pi i}\right)^d\right] \in
H^{2d}(X;\Q)\] is the $d$th part of the Chern character of $E$.

For a holomorphic bundle $E$ on a complex manifold $X$ there is a
preferred class of connections coming from smooth hermitian metrics on
the bundle. These define curvature forms of type $(1,1)$. Now an
$\End(E)$-valued $(1,1)$-form can also be thought of as a
$(T^\dualstar\otimes\End(E))$-valued $(0,1)$-form, where $T^\dualstar$
is the holomorphic cotangent bundle of $X$. After this identification
we are free to use {\em more complicated operations} to combine the
curvature with itself (as well as with the curvature of the
holomorphic tangent bundle and a holomorphic symplectic form, if
available), because the curvature now has three tensorial ``indices''
rather than the original two. The different possible combinations,
which replace the invariant polynomials used above, are in fact
parametrised by Jacobi diagrams such as those defining $\A$.


\subsection{Curvature of a holomorphic bundle}

In order to fix the notation, let us recall the basics of complex
differential geometry.  If $X$ is a complex manifold then one may
decompose the complexified tangent bundle into holomorphic and
anti-holomorphic parts: $T_\R X \otimes \C \cong T\oplus \bar T$.  The
exterior differential likewise splits as $d=\partial + \bar \partial$
and the complexified de Rham complex $(\Omega^\gradestar (X;\C), d)$
may be refined to obtain the Dolbeault complex $(\Omega^{\gradestar
,\gradestar }(X;\C), \bar\partial)$, with cohomology $H^{\gradestar
,\gradestar }_{\bar\partial}(X;\C)$. (If $X$ has a K\"ahler metric
then these Dolbeault cohomology groups may be identified with
subspaces of the complexified de Rham cohomology of $X$, via the Hodge
decomposition, but we will not usually assume $X$ is K\"ahler in what
follows.)

Unlike in the smooth case, if $E$ is a holomorphic vector bundle on
$X$, then there is a canonical operator $\bar \partial$ on the spaces
of smooth $E$-valued forms $\Omega^{p,q}(X;E)$. To define it, write a
form locally in terms of a basis of holomorphic sections of $E$, and
apply the usual Dolbeault $\bar \partial$ operator to the smooth-form
coordinates; one obtains a complex with cohomology
$H^{\gradestar,\gradestar}_{\bar\partial}(X;E)$.

A smooth connection on a holomorphic bundle $E$, thought of as a
covariant derivative $\nabla: \Omega^0(X;E) \ra \Omega^1(X;E)$,
splits into pieces of type $(1,0)$ and type $(0,1)$. It is said to be
{\em compatible} with the holomorphic structure on $E$ if its $(0,1)$
part equals the canonical $\bar \partial$ operator of $E$. One may
then write $\nabla = \bar \partial + \nabla^{1,0}$, the last term
being a\/ {\em connection of type $(1,0)$}, which satisfies a version
of the usual Leibniz rule in which $\partial$ replaces $d$. The
resulting curvature two-form $F \in \Omega^2(X;\End(E))$ has no part
of type $(0,2)$, because $F = (\bar\partial + \nabla^{1,0})^2$ and
$\bar \partial ^2=0$. In local coordinates, if one writes the
covariant derivative operator $\nabla$ as $d+A$ for some $1$-form $A
\in \Omega^1(X;\End(E))$, this compatibility amounts to saying that
$A$ is of type $(1,0)$.

If $E$ has a smooth hermitian metric $h$ then we may further require
that $\nabla$ is compatible with $h$ by imposing that for all sections
$s,t \in \Omega^0(X;E)$,
\[ d h(s,t) = h(\nabla s,t) + h(s, \nabla t).\]
Computing $d$ of this formula using a basis of local
covariant-constant sections shows that the curvature $F=\nabla^2$ is
of type $(1,1)$ (and in fact purely imaginary). Therefore
$(\nabla^{1,0})^2=0$ and we can write the operator $F$ as
$\nabla^{1,0} \bar \partial + \bar \partial \nabla^{1,0}$ or even as
$\nabla \bar \partial + \bar \partial \nabla$.  This second form will
be used below.  Varying the hermitian form alters the form $F$ by a
$\bar\partial$-coboundary.

If two bundles $E_1,E_2$ have connections, then there is an induced
connection on $E_1 \otimes E_2$ given by the Leibniz rule, and the
resulting curvature is \[ F_{E_1 \otimes E_2} = F_{E_1} \otimes \id +
\id \otimes F_{E_2}.\] Similarly, a connection on a bundle $E$ induces
one on its dual $E^\dualstar $ by the formula
\[ \langle\nabla\phi, s\rangle + \langle \phi, \nabla s\rangle =d
\langle \phi, s \rangle, \] where $s$ is a section of $E$, $\phi$ is a
section of $E^\dualstar$ and the brackets indicate the contractions to
complex valued forms on $X$. It is useful to think in terms of
operators on the space of sections of $E$ and write $F_{E^\dualstar
}\phi = -\phi \circ F_E$.

The {\em Bianchi identity} is often written $\nabla F=0$. The operator
$\nabla = \nabla_{\End(E)}$ is the covariant derivative on sections of
the bundle $\End(E)$ induced by the original connection $\nabla$ on
$E$. As an operator on $\Omega^0(X;E)$, $\nabla_{\End(E)} F=\nabla_E
\circ F - F \circ \nabla_E$, so its vanishing amounts to nothing more
than the fact that the operators $F=\nabla_E^2$ and $\nabla_E$
commute.  In the holomorphic context, the $(1,0)$ part of the identity
becomes the equation $\bar \partial F=0$.


\subsection{Complex manifolds and the Jacobi identity}

Kapranov discovered that the curvature of a holomorphic bundle on a
complex manifold satisfies a kind of Jacobi identity. This fact (which
has nothing to do with hyperk\"ahler or complex symplectic geometry)
is absolutely basic to Rozansky-Witten theory.

Suppose $E$ is a holomorphic bundle on $X$, with associated Dolbeault
operator $\bar \partial_E$. Pick a smooth hermitian metric on $E$
with associated connection $\nabla_E$ and curvature form $F_E\in
\Omega^{1,1}(X,\End(E))$.  Do the same for the holomorphic tangent
bundle $T$.  We will from now on drop the redundant ``$X$'' from
notation such as $\Omega^\gradestar (X;E)$.

We want to think of the curvature as living in a slightly different
space. Let $\Theta$ denote any identification of the form
$\Omega^{p,q}(-) \cong \Omega^{0,q}(\bigwedge^p(T^\dualstar ) \otimes
-)$. Here we think of the right hand side as a subspace of
$\Omega^{0,q}((T^\dualstar )^{\otimes p} \otimes -)$, and explicitly
(this will affect signs in an inevitably messy way) set $\Theta (d\bar
z^I \wedge dz^J \otimes s) = d\bar z^I \otimes dz^J \otimes s$.
Define $R_E =\Theta F_E \in \Omega^{0,1}(T^\dualstar \otimes
\End(E))$; this form will also be referred to as the curvature.
Since $F_E$ is $\bar\partial$-closed, so is $R_E$, as the appropriate
$\bar\partial$ operators commute with $\Theta$.  Define $R_T$
similarly.

Kapranov's result is that a certain three-term quadratic relation in
the tensors $R_E, R_T$ is a $\bar\partial$-coboundary. At the level of
cohomology it will become the {\em STU relation} of Vassiliev theory,
and in the special case $E=T$ the {\em IHX relation}.  Define three
elements of $\Omega^{0,2}(T^\dualstar \otimes T^\dualstar \otimes
\End(E))$ called $R_E \circ_S R_T, R_E \circ_T R_E, R_E \circ_U R_E$
by taking the appropriate wedge products of $1$-forms and contracting
indices according to the three graphs shown below.
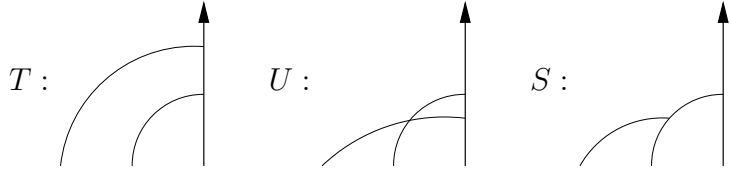
\[ T:\vpic{curvT} \qquad U:\vpic{curvU} \qquad S: \vpic{curvS}\]
Explicitly, applying these elements to sections $t_1$, $t_2$, and $e$
one gets elements of $\Omega^{0,2}(E)$ which may be written
$R_E(R_T(t_1, t_2), e), R_E(t_1, R_E(t_2, e)), , R_E(t_2, R_E(t_1,
e))$.

\begin{lemma}[STU relation]
If $E$ is a holomorphic bundle over the complex manifold $X$, then in
$\Omega^{0,2}(T^\dualstar \otimes T^\dualstar \otimes \End(E))$ we
have the coboundary formula: \[ R_E \circ_T R_E + R_E \circ_U R_E +
R_E \circ_S R_T=-\bar\partial (\Theta \nabla R_E).\]
\end{lemma}
\begin{proof}
Via the Leibniz formula we obtain the operator identity \[
F_{T^\dualstar \otimes \End(E)} = F_{T^\dualstar \otimes E^\dualstar
\otimes E} = F_{T^\dualstar } \otimes \id \otimes \id + \id \otimes
F_{E^\dualstar } \otimes \id + \id \otimes \id \otimes F_E, \] so that
composing with $R_E$ and evaluating on sections $t,e$ of $T,E$ we have
in $\Omega^{1,2}(E)$ the identity \[ (F_{T^\dualstar \otimes
\End(E)}R_E)(t,e) = - R_E(F_Tt, e) - R_E(t, F_Ee) + F_E (R_E(t,e)).  \]
(The signs come from the curvature of the dual bundle; switching the
order of $2$-form and $1$-form does not give signs.) Now applying
$\Theta$ (carefully) to obtain an identity in
$\Omega^{0,2}(T^\dualstar \otimes T^\dualstar \otimes
\End(E))$ gives \[ -\Theta(F_{T^\dualstar \otimes \End(E)}R_E) = R_E
\circ_S R_T + R_E \circ_T R_E + R_E \circ_U R_E.  \] The result now
follows on rewriting the left-hand side using $F R_E=(\bar\partial
\nabla + \nabla \bar\partial) R_E =\bar\partial (\nabla R_E)$ (because
$R_E$ is $\bar\partial$-closed) and the fact that $\Theta$ commutes
with $\bar\partial$.
\end{proof}

Just as important from the point of view of constructing weight
systems is the {\em symmetry} of the curvature form $R_T$ of the
tangent bundle. In fact there are two separate symmetries: the first
comes from considering the torsion of the connection on $T$, while the
second appears in the presence of a holomorphic symplectic form, and
will be studied in the next section. Kapranov assumes in his paper
that the hermitian metric on $X$ is K\"ahler, so that the torsion of
$\nabla_T$ vanishes (this is one definition of a K\"ahler metric, in
fact). But the next proposition shows that vanishing of the torsion is
unnecessary; one no longer has an exact symmetry, but symmetry modulo
coboundaries, which is still perfectly acceptable to us.

If $\nabla$ is a smooth connection on the real tangent bundle $T_{\R}$
of a smooth manifold, then the {\em torsion} is a $2$-form with values
in $T_{\R}$ given by the formula
\[\tau(t_1,t_2)=\nabla_{t_1} t_2 -\nabla_{t_2} t_1 -[t_1,t_2].\]
For a complex manifold with a smooth hermitian connection $\nabla$ on
its holomorphic tangent bundle $T$, we can tensor over $\R$ with $\C$
to obtain a connection all of $\T_{\C}=T\oplus \overline{T}$, and use
the same formula to define the torsion $\tau \in
\Omega^2(T_{\C})$. The part $\tau^{1,0}$ with values in $T$ turns out
to be of type $(2,0)$.

\begin{prop}[Partial symmetry]
The curvature form $R_T$ is symmetric in its two inputs, up to a $\bar
\partial$-coboundary. Specifically, \[ R_T -\sigma \circ R_T
=\bar\partial ( \Theta\tau^{1,0}) \] where $\Theta\tau^{1,0} \in
\Omega^{0,0}(T^\dualstar \otimes T^\dualstar \otimes T)$ is a version
of the torsion and $\sigma$ is the permutation of the two $T^\dualstar
$ factors.
\end{prop}
\begin{proof}
For arbitrary smooth sections $t_1, t_2, t_3$ of $T_{\C}$ we have the
straightforward identity
\[ \sum F(t_1,t_2)t_3 = \sum d_{t_1} \tau(t_2,t_3) + \sum \tau
(t_1, [t_2,t_3]) = (\nabla \tau)(t_1, t_2, t_3),\] all sums being over
cyclic permutations of the three vector fields. (In the Levi-Civita
case, the vanishing of the RHS implies one of the symmetries of the
Riemann curvature.) Now assume $t_1, t_3$ are of type $(1,0)$ while
$t_2$ is of type $(0,1)$, and look at the type $(1,0)$ part of this
equation:
\[ F(t_1,t_2)t_3 + F(t_2,t_3)t_1 = (\nabla \tau^{1,0})(t_1, t_2, t_3)
= (\bar\partial \tau^{1,0})(t_1, t_2, t_3).\] Applying $\Theta$ we can
have an identity in $\Omega^{0,1}(T^\dualstar \otimes T^\dualstar )$
which when evaluated on $t_2,t_1,t_3$ says that
\[ R_T(t_2)(t_1,t_3) - R_T(t_2)(t_3,t_1) = \Theta(\bar\partial
\tau^{1,0})(t_2)(t_1,t_3),\]
as required.
\end{proof}

\begin{remk}
The exterior product of forms followed by contraction with $R_T$
defines a degree-one bilinear product on the Dolbeault complex
$\Omega^{0,*}(T)$. This operation satisfies the graded Jacobi identity
up to a coboundary, and in the K\"ahler case it is exactly symmetric,
making it an ``odd Lie bracket up to homotopy''. Kapranov shows that
together with higher-order derivatives of the curvature, it makes the
Dolbeault complex $\Omega^{0,*}(T)$ into an {\em
$L_\infty$-algebra}. In the non-K\"ahler case, the above lemma
suggests that there is an even weaker kind of infinity-structure in
which there are also higher homotopies (controlled by derivatives of
the torsion) arising from non-commutativity of the bracket. Such
structures are beautiful and interesting, but we will not need them in
this paper.
\end{remk}


\subsection{Complex symplectic manifolds}

As we have seen, the curvature of a holomorphic vector bundle has a
kind of intrinsic Jacobi identity property. To construct weight
systems we also need a {\em metric} of some kind, and in keeping with
the ``switch of statistics'' that has replaced a skew Lie bracket by a
symmetric curvature tensor, we seek a {\em skew} rather than symmetric
non-degenerate form.

A {\em complex symplectic manifold} is an (even-dimensional) complex
manifold $X^{2n}$ together with a non-degenerate holomorphic two-form
$\omega\in \Omega^0(\bigwedge^2 T^\dualstar )$. The non-degeneracy
implies that $\omega$ defines an isomorphism of holomorphic bundles $T
\cong T^\dualstar $. (An obvious topological obstruction to existence
is therefore the vanishing of the odd rational Chern classes of $X$.)

Using this isomorphism, we can convert the curvature $R_T \in
\Omega^{0,1}(T^\dualstar \otimes
\End(T))$ into a form $C_T \in \Omega^{0,1}(T^\dualstar \otimes
T^\dualstar \otimes T^\dualstar )$:
\[ C_T(t_1,t_2,t_3) = \omega(R_T(t_1,t_2), t_3). \]
Since $\omega$ is holomorphic, $C_T$ too is $\bar\partial$-closed.

\begin{lemma}[Full symmetry of curvature]
The curvature form $C_T\in \Omega^{0,1}(T^\dualstar \otimes
T^\dualstar \otimes T^\dualstar )$ of a complex symplectic
manifold is symmetric in its three factors, up to
$\bar\partial$-coboundaries.
\end{lemma}
\begin{proof}
We already have such symmetry in the first two factors. To show
symmetry in the second and third, consider $F_{T^\dualstar \otimes
T^\dualstar }\omega$.  By the Leibniz rule and the rule for curvature
of dual bundles we can write \[ F_{T^\dualstar \otimes T^\dualstar
}\omega = -\omega \circ (F_T \otimes \id_T) -\omega \circ (\id_T
\otimes F_T).  \] Applying $\Theta$ and rewriting this identity in
terms of elements of $\Omega^{0,1}(T^\dualstar \otimes T^\dualstar
\otimes T^\dualstar )$, gives \[ \Theta(\bar\partial (\nabla \omega))
= -C_T +C_T \circ \sigma_{23}, \] where $\sigma_{23}$ is the
permutation of the last two inputs. The left-hand side is the
coboundary $\bar \partial (\Theta (\nabla\omega))$ and so the symmetry
is proved. (Note that $\nabla$ came from an arbitrary choice of
hermitian metric on $T$; there is no reason why $\nabla \omega$ should
be zero.)
\end{proof}


\subsection{Rozansky-Witten weight systems}

With the above preliminaries completed, we can now describe briefly
the construction of weight systems on the space $\A$.

\begin{thm}
If $X$ is a complex symplectic manifold and $E$ is a holomorphic
vector bundle on $X$, then there is a weight system
\[ RW_{X,E}:\A \ra H_{\bar\partial}^{0,\gradestar }(X).\]
taking values in the Dolbeault cohomology of $X$.
\end{thm}
\begin{proof}
If $\Gamma$ is a $2v$-vertex closed trivalent graph with {\em
ordered vertices} and {\em oriented edges}, then one can obtain a
form in $\Omega^{0,2v}(X)$ by a procedure like that of
Section~\ref{Section:LAweightsystems}: wedge/tensor one copy of
$C_T$ for each vertex of $\Gamma$, and contract tensorially with
one copy of $\omega^{-1}$ (that is, $\omega$ converted into a
holomorphic section of $T \otimes T$) for each edge of $\Gamma$.
Since all the elements in the construction are
$\bar\partial$-closed, so is the result, by the Leibniz formula.
There is clearly a choice of how one attaches $C_T$ at a vertex
--- a choice of correspondence between the three legs and the
three copies of $T^\dualstar $ --- but differences alter the
resulting form by a coboundary, because of the symmetry property
and the Leibniz formula. The choice of hermitian metric used to
define $C_T$ similarly only affects the result by a coboundary, so
that the result is a well-defined element of
$H_{\bar\partial}^{0,2v}(X)$.  This basic construction clearly
generalises to the case where the graph has an oriented Wilson
loop: the form $R_E$ is inserted at the vertices on the loop,
which are canonically oriented.

Reversing the orientation of any edge or swapping the order of two
vertices negates this element, because $\omega^{-1}$ and the cup
product of $1$-forms are skew. (The following example may help: if
$\alpha, \beta$ are $1$-forms with values in vector space $V,W$,
then $\alpha \wedge \beta = - \sigma(\beta \wedge \alpha)$, where
$\sigma$ is the usual permutation. In particular $\alpha \wedge
\alpha = -\sigma(\alpha \wedge \alpha)$; ``$\alpha$ anticommutes
with itself''.) Therefore the map is really well-defined on {\em
oriented} graphs, where an orientation is an ordering of the
vertices and an orientation of the edges, considered up to an even
number of transpositions and reversals. The remarkable fact is
that this notion of orientation is canonically isomorphic to the
standard convention (from Section~\ref{Section:LAweightsystems})
on Jacobi diagrams, in which each vertex has a cyclic ordering of
its legs.

To see this, let $V$ and $E$ be, respectively, the sets of vertices
and edges of $\Gamma$. Let $F$ be the set of all flags (half-edges) of
$\Gamma$, and for each vertex $v$ and for each edge $e$ let $F_v$ and
$F_e$ be the obvious two- and three-element sets of incident flags.
For any set $S$, use the notation $\Det(S)$ for the top exterior power
$\Det(\R^S)$, so that orienting a vertex or an edge in the usual sense
amounts to orienting the appropriate $1$-dimensional vector space
$\Det(F_v)$ or $\Det(F_e)$. Orientations of graphs under the two
different conventions are measured by the spaces $\Det(V) \otimes
\bigotimes_e \Det(F_e)$ and $\bigotimes_v \Det(F_v)$.

The isomorphism now follows by combining three simple natural
(equivariant) isomorphisms: (i) $\Det(F) \cong \Det(V) \otimes
\bigotimes_v \Det(F_v)$; (ii) $\Det(F) \cong \bigotimes_e \Det(F_e)$;
and (iii) $\Det^2$ is canonically trivial.  The first two isomorphisms
come from concatenating triples of flags according to the vertex
order, or pairs of flags according to an arbitrary (irrelevant) edge
order.

The fact that the construction respects the IHX and STU relations now
follows from the earlier proposition about the Jacobi identity for the
curvature.
\end{proof}

\begin{remk}
This last check is actually quite nasty, because the two different
orientation conventions we are considering do not agree {\em locally},
and the equivalence between them is not so straightforward even {\em
globally}. The categorical approach we adopt in the second half of the
paper has a technical advantage in that it matches the correct
orientation conventions {\em locally}, bypassing this annoying
problem.
\end{remk}

\subsection{Examples}

In many ways the best examples of complex symplectic manifolds are the
{\em hyperk\"ahler manifolds}, which were the subject of Rozansky and
Witten's original work. A hyperk\"ahler manifold is a real
$4n$-dimensional manifold with a Riemannian metric of holonomy
$\Sp(n)$. Because this group is contained in $GL(n,\HHH)$, one can
introduce three parallel (which implies integrable) almost complex
structures $I$, $J$, and $K$ satisfying the usual quaternionic
relation $IJK=-1$. Any imaginary unit quaternion $q$ now defines a
complex structure (for which the metric is K\"ahler, with K\"ahler
form $\omega_q$) and which possesses a holomorphic symplectic form:
one only needs to check for example that the complex two-form
$\omega=\omega_J+i\omega_K$ is $I$-holomorphic.

There is a partial converse: a {\em compact} complex symplectic
manifold which is K\"ahler has a hyperk\"ahler metric, by Yau's
solution of the Calabi conjecture (Beauville \cite{Beauville}). (This
is a hard analytical existence theorem, and there is no known simple
formula for the metric.) Kapranov's approach is therefore only really
more general than Rozansky and Witten's if we are prepared to consider
complex symplectic manifolds which are non-compact, non-K\"ahler, or
both.  There are a few compact non-K\"ahler examples due to Beauville
and Guan \cite{Guan}, but there are plenty of non-compact
hyperk\"ahler manifolds coming from complex Lie group coadjoint
orbits, geometric moduli spaces, etc. (see Hitchin \cite{H}).

From the point of view of Vassiliev invariants, the {\em compact} case
is (at least initially) the most interesting, because for a compact
complex symplectic manifold $X$ of real dimension $4n$, one can obtain
{\em scalar}-valued weight systems of degree $2n$. To do this,
integrate the invariants lying in $H_{\bar\partial}^{0,2n}(X)$ against
the holomorphic volume form $\omega^n \in H_{\bar\partial}^{2n,0}(X)$.
Further, in the {\em hyperk\"ahler} case, Sawon \cite{Sawonthesis}
used the interplay between the Riemannian and hermitian constructions
to show that these numbers are invariant under deformations of the
hyperk\"ahler metric and of the complex structure on $X$. He also
performed some explicit calculations.

The current list of known compact hyperk\"ahler manifolds is not very
long. In dimension four, the K3 surface and $4$-torus are the only
examples. Each of these generates, via its Hilbert schemes of points
(desingularised versions of its symmetric products), an infinite
family of further examples. These are all {\em irreducible}, having
holonomy not contained in a proper subgroup of $\Sp(n)$, and in
particular not being products of lower dimensional hyperk\"ahler
manifolds. The only other known irreducible examples were both
constructed by O'Grady \cite{OG1, OG2}.

The relative paucity of examples --- two countable families and some
exceptions --- might therefore seem to undermine the scope of the
Rozansky-Witten weight systems. But in fact if one looks to Lie
algebras one finds exactly the same situation --- the two series of
types $A$ and $BCD$, and a few exceptions! In this sense there are
``at least as many'' Rozansky-Witten weight systems as Lie algebra
ones.  An obviously important issue is whether the Rozansky-Witten
weight systems are really {\em new}, lying outside the span of the Lie
algebra ones or not. Because of the difficulties in explicit
calculation, we don't yet know the answer to this.

\pg


\section{Lie algebra weight systems revisited}
\label{Section:LAweightsystemsAgain}
 In this section we describe an
alternative category-theoretic approach to the construction of
weight systems from metric Lie algebras. It was introduced by
Vogel \cite{Vogel} and Vaintrob \cite{Vaintrob}, whose original
motivation was to handle the weight systems arising from metric
Lie superalgebras.

For such an algebra, the tensors $f$ and $b$ used in
Section~\ref{Subsection:LAweightsystems} have both skew and
symmetric parts, leading to incompatibility with the standard
orientation convention for Jacobi diagrams. The problem can be
fixed by picking a direction in the plane and representing Jacobi
diagrams always as Morsified planar graphs, rather than as
abstract graphs. The approach leads inevitably to the idea of
constructing weight systems from metric Lie algebras in any
category for which the notion makes sense, and not just in the
category of (super-)vector spaces. We will justify all this
abstract nonsense later in the paper by constructing interesting
{\em examples} of such categories and Lie algebras.


\subsection{Symmetric tensor categories.}
\label{Subsection:SymmetricTensorCats} Here we will recall the
standard definitions of symmetric tensor categories. For more
detail see Bakalov and Kirillov \cite{BK}, Chari and Pressley
\cite{CP}, or Kassel \cite{Kassel}.

A category $\scrC$ is a {\em tensor} (or {\em monoidal}) category
if it comes with a functor $\otimes: \scrC \times \scrC \ra \scrC$
whose associativity is implemented by a natural isomorphism $\Phi
: \otimes \circ (\otimes \times \id) \ra \otimes \circ (\id \times
\otimes)$ satisfying the pentagon identity, and has a unit object
$1$ for tensor product, again with appropriate natural
isomorphisms. We will for the moment ignore all these isomorphisms
notationally, pretending that $\scrC$ is strictly associative,
i.e. that these isomorphisms are equalities.

A symmetric tensor category is defined as follows.  Let $\sigma$
be the standard flip functor $\scrC \times \scrC \ra \scrC \times
\scrC$; $A\otimes B\mapsto B\otimes A$. The tensor category
$\scrC$ is {\em symmetric} if there is a natural isomorphism
$\tau\colon \otimes \ra \otimes \circ \sigma$ --- giving an
isomorphism $\tau_{A,B}\colon A\otimes B\to B\otimes A$ for $A,B$
objects of $\scrC$ --- satisfying $\tau_{B,A}\circ\tau_{A,B}=\id$
and satisfying the hexagon relation
\[  \tau_{A, B \otimes C} =(\id_B \otimes \tau_{A,C}) \circ (\tau_{A,B}
\otimes \id_C). \] The hexagon would be more visible if we hadn't
dropped the associators from the notation.  The natural
isomorphism  $\tau$ is sometimes called the {\em symmetry}. The
standard example to keep in mind here is the category of {\em
super-vector spaces}, which is symmetric but in a non-trivial way;
the isomorphism $\tau$ will handle all the signs for us.

The notion of {\em duality} between objects in a tensor category
$\scrC$ is a little tricky. The basic definitions are abstracted
from properties of finite-dimensional vector spaces, but a little
more is required in order to control double duals properly. An
object $A^\dualstar $ is a {\em right dual} of an object $A$ if
there is a\/ {\em right evaluation\/} $\epsilon_A\colon
A^\dualstar \otimes A\ra 1$ and a\/ {\em right co-evaluation\/}
$\iota_A\colon 1 \ra A \otimes A^\dualstar $ which satisfy
\begin{gather*}
  (\id_A \otimes \epsilon_A) \circ (\iota_A \otimes\id_A) = \id_A\\
  (\epsilon_A \otimes \id_{A^\dualstar }) \circ
  (\id_{A^\dualstar } \otimes \iota_A) = \id_{A^\dualstar }.
\end{gather*}
Such an object is unique up to a canonical isomorphism. We can
similarly define a left dual ${}^\dualstar\!\! A$ with structural
maps $\iota_A'\colon \one \ra {{}^\dualstar\!\! A} \otimes A$ and
$\epsilon_A'\colon A\otimes {{}^\dualstar\!\! A}\ra \one$.  A {\em
rigid} tensor category is one in which all objects have left and
right duals. This is enough to permit the construction of traces
$\Mor(A,A) \ra \Mor(\one,\one)$ on endomorphisms of any object,
the construction of adjoints of morphisms, and identifications
such as $\Mor(A\otimes B, C)\cong\Mor(A, C\otimes B^\dualstar )$.

Most of the categories we will use in this paper will be at least {\em
additive} (and probably {\em $\C$-linear}), having abelian groups (or
complex vector spaces) for morphism sets, bilinear composition, a
direct sum operation $\oplus: \scrC \times \scrC \ra \scrC$ and a
zero object $0$.


\subsection{Penrose calculus}

An important tool is Penrose's diagrammatic representation of the
structure of a tensor category by planar pictures. A tensor product of
objects is represented by a collection of labelled dots on a
horizontal level; a morphism from one such to another is represented
by drawing, inside a horizontal strip whose top and bottom edges are
labelled appropriately, a box, labelled with the name of the morphism,
and connected by strings from its top and bottom edges to the object
dots. Composition of morphisms is represented by concatenation of
diagrams moving up the page; tensor product of morphisms by horizontal
juxtaposition. (If we were not assuming strict associativity then
bracketings of objects and explicit associator morphisms would also be
required.)

Special structural morphisms in the category are represented using
special pictures as a short-hand for labelled boxes. The identity
morphism on an object is always represented by a vertical arc
labelled with that object, and the other possible structural
morphisms are pictured below. The point of using these particular
pictures is of course that the rather complicated algebraic {\em
relations} satisfied by the structural morphisms now correspond to
natural topological identities.

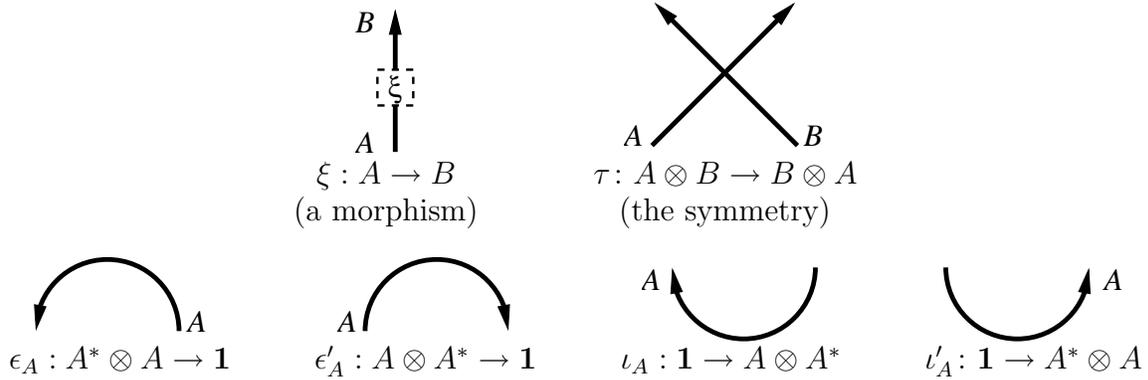
\[\begin{array}{ccc}\vpic{btrxi}&\qquad&\vpic{tau}\\
\xi: A \ra B&&\tau\colon A\otimes B\to B\otimes A\\
\text{(a morphism)}&&\text{(the symmetry)}
\end{array}\]
\[\begin{array}{ccccccc}
\vpic{capepsL}&\quad&\vpic{capepsR}&\quad&
\vpic{cupepsL}&\quad&\vpic{cupepsR}\\
\epsilon_A:  A^\dualstar \otimes A\ra\one&& \epsilon_A':  A
\otimes A^\dualstar  \ra \one&& \iota_A: \one\ra A\otimes
A^\dualstar && \iota_A'\colon  \one \ra A^\dualstar \otimes A
\end{array}
\]


\subsection{Lie algebras and modules}
Here we take the usual definitions of a Lie algebra and a Lie
algebra module and abstract them from the category of vector
spaces to an arbitrary additive symmetric tensor category.

 Let $\scrC$ be an additive
symmetric tensor category.  A {\em Lie algebra} in $\scrC$ is an
object $L$ equipped with a {\em bracket} morphism $\alpha\colon
L\otimes L\to L$ which is skew-symmetric and satisfies the Jacobi
identity: %
\begin{gather*}
  \alpha+\alpha\circ \tau=0;\\
  \alpha\circ(\alpha\otimes \id)+ \alpha\circ(\alpha\otimes
\id)\circ
  \tau_{123} +\alpha\circ(\alpha\otimes \id)\circ
  \tau_{321}=0;        
\end{gather*} where $\tau_{123}$ and $\tau_{321}$ denote the
actions on $L^{\otimes 3}$ of the the three-cycles in the
symmetric group $S_3$.  Note that addition of morphisms makes
sense because $\scrC$ is additive.

A (right) {\em module} over such a Lie algebra is an object $M$
together with an {\em action} morphism $\alpha_M\colon M\otimes L
\to M$ satisfying the identity
\[   \alpha_M\circ(\id\otimes \alpha)= \alpha_M\circ(\alpha_M\otimes
\id)-\alpha_M\circ(\alpha_M\otimes \id)\circ (\id\otimes \tau).\]
Pictorially, the bracket and action are represented by the following
diagrams, in a way that turns the above identities into the
antisymmetry, IHX and STU relations.
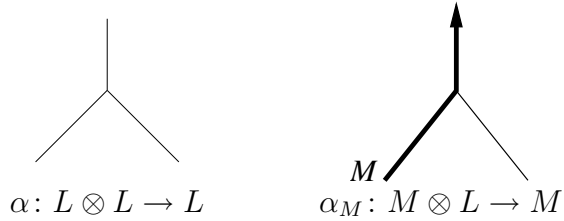
\[\begin{array}{ccc}
\vpic{alpha}&\qquad&\vpic{alphaM}\\
\alpha\colon  L\otimes L\to L&&\alpha_M \colon  M\otimes L\to M
\end{array}\]
Note that any right module $M$ can be given a natural left module
structure $\check{\alpha}_M\colon L\otimes M\to M$ by
$\check{\alpha}_M=-\alpha_M\circ\tau_{L,M}$.

 An {\em $L$-module
morphism} is a morphism $\xi: M\to N$ between $L$-modules such
that $\xi\circ \alpha_M=\alpha_N\circ (\xi\otimes \id)$.
Pictorially this is shown below.  The collection of $L$--modules
and $L$--morphisms form a category $\modL$.
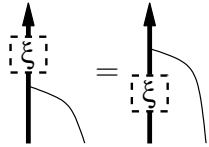
\[ \vpic{xiL} = \vpic{xiR}\]

The {\em tensor product} of two $L$-modules is an $L$-module under
a Leibniz rule such as $\alpha_{M\otimes N}=(\alpha_M\otimes
\id)\circ(\id \otimes \tau)+ \id\otimes\alpha_N$, and therefore
$\modL$ is a tensor-category.  The action on a tensor product is
defined and notated as shown below. The crossings have been drawn
in a slightly non-Morse way here, but we hope that the meaning is
clear: they are $\tau$ morphisms forming an essential part of the
correct definition of the action on tensor products.
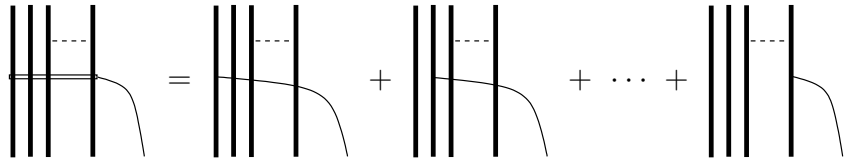
\[  \vpic{Leibniz} ~=~ \vpic{Leibniz1} ~+ ~\vpic{Leibniz2}~ + ~\cdots~ +~
\vpic{Leibnizn}. \]

 A {\em metric Lie algebra} is a Lie algebra
equipped with an abstracted version of a non-degenerate symmetric
invariant bilinear form. Thus, it comes with a {\em metric}
morphism $\beta: L\otimes L \ra \one$ and a {\em Casimir} $\gamma:
\one \ra L\otimes L$, each an $L$-module morphism satisfying
non-degeneracy and symmetry axioms:
\begin{gather*}
  (\id\otimes \beta)\circ(\gamma\otimes \id) = \id =
  (\beta \otimes \id)\circ(\id\otimes \gamma);\\
  \beta=\beta\circ\tau; \qquad \gamma = \gamma \circ \tau.
\end{gather*}

In pictures, cup and cap denote these morphisms
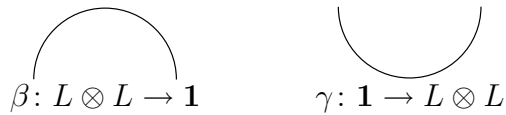
\[  \begin{array}{ccc}
\vpic{beta}&\qquad&\vpic{gamma}\\
\beta\colon  L\otimes L\ra \one&&\gamma\colon  \one\ra  L\otimes L
\end{array}
\]

When $\scrC$ is a rigid category, the dual of a module may be made a
module by forcing the evaluation and coevaluation maps to be module
maps. This is better defined by a picture than by a formula:
\[ \vpic{dual1}~ = ~-\vpic{dual2}.\]

\subsection{Weight systems}

With this framework set up, we can state the theorem which will
underlie our later explicit construction of the Rozansky-Witten weight
systems:

\begin{thm}[\cite{Vaintrob,Vogel}]
Let $\scrC$ be a rigid, additive, symmetric tensor category, $L$ a
metric Lie algebra in $\scrC$, and $M$ a dualizable module over
$L$. Then there is a weight system
\[ w_{L,M}: \A \ra \Mor(\one,\one).\]
\end{thm}

\begin{proof}
Given any Jacobi diagram in $\A$, first draw it in the plane in a way
compatible with its orientation. Morsify it so that the critical
points and trivalent vertices lie at different levels, and so that the
whole diagram is built from the generating morphisms we gave earlier,
together with the Lie bracket and module action. Now compose the
corresponding morphisms in $\scrC$. The proof of independence of the
Morse and planar structures is the usual Reidemeister-move type
argument, for which we refer to Vaintrob \cite{Vaintrob}.
\end{proof}



\pg

\section{Sheaves and derived categories}
\label{Section:SheavesAndDerivedCategories}

Our main goal in this paper is to reinterpret the Rozansky-Witten
weight systems in the context of the category-theoretic framework
described above. The basic construction will be to associate to any
complex manifold $X$ a symmetric tensor category $D(X)$ and a Lie
algebra object $L$ in $D(X)$.

We begin in this section with a quick general explanation of the
salient points about derived categories and about sheaves.  Useful
references for derived categories are Gelfand and Manin \cite{GM} and
Richard Thomas \cite{Th}. For sheaves see Hartshorne \cite{Ha} and
Kashiwara-Schapira \cite{KS}. 

\subsection{Derived categories}

Let $\scrC$ be an {\em abelian category}; recall that this is an
additive category in which every morphism has a kernel and a cokernel,
and the two possible definitions of ``image'' (cokernel of kernel, or
kernel of cokernel) agree. The standard example is the category of all
(right, say) modules over a ring $R$, and in practice one may treat
any abelian category as being of this form. From $\scrC$ we can form
the category $\Ch(\scrC)$ of {\em chain complexes} of objects of
$\scrC$.

In homological algebra, one works primarily at the level of chain
complexes, because taking homology groups prematurely can destroy some
of the information they contain. For example, the homology-cohomology
universal coefficient theorem shows that the operation of replacing a
complex by its homology does not commute with the operation of taking
the dual. When working in $\Ch(\scrC)$, it is clearly reasonable to
identify chain-homotopic maps and thereby to pass to a quotient {\em
homotopy category} $\Hty(\scrC)$, whose morphisms are the homotopy
classes of maps between complexes.

But it is more sensible to regard in addition any {\em
quasi-isomorphism} -- a map between complexes which induces
isomorphisms on homology -- as an isomorphism. Although chain homotopy
equivalences are certainly quasi-isomorphisms, the converse is not
true; there may remain in $\Hty(\scrC)$ quasi-isomorphisms without
inverses. This can cause problems: for example, we often want to view
a module as ``equivalent'' to any of its projective resolutions
(quasi-isomorphic complexes of projective modules), but such
resolutions need not actually be homotopy-equivalent to the original
module.

The {\em derived category} $D(\scrC)$ is defined by formally
inverting these inside $\Hty(\scrC)$: one introduces a calculus of
fractions $f/g$ (for $f$ any morphism and $g$ a quasi-isomorphism)
essentially identical to the Ore localisation for non-commutative
rings. Explicitly, any morphism in $D(\scrC)$ between the complexes
$A^\bl$ and $B^\bl$ may be represented by a diagram of each of the
forms
\[ A^\bl \stackrel{f}{\ra} C^\bl \stackrel{g}{\leftarrow}
B^\bl\qquad \text{and} \qquad A^\bl \stackrel{g}{\leftarrow} C^\bl
\stackrel{f}{\ra} B^\bl, \] for some other complex $C^\bl$.

Any functor defined on $\Ch(\scrC)$ which takes quasi-isomorphisms to
isomorphisms -- the abelian-group-valued homology functors $h^i:
\Ch(\scrC) \ra \Ab$ being the obvious examples -- therefore factors
through $D(\scrC)$, and in fact this universal property characterises
$D(\scrC)$.

Note that the objects of the derived category are the same as those of
$\Ch(\scrC)$, and that objects of the original category $\scrC$ may be
identified with chain complexes whose only non-zero term lies in
degree $0$, so that there is an ``inclusion'' functor $\scrC \ra
D(\scrC)$.

In $D(\scrC)$ there are {\em shift functors} written $[n]\colon
D(\scrC) \ra D(\scrC)$, for $n\in \Z$. The functor $[n]$ acts on a
complex $A^\bl$ by shifting it $n$ places to the left, so that
$A^\bl[n]^i = A^{i+n}$ and the differential is $d[n]^i=(-1)^n
d^{i+n}$. It acts on chain maps by shifting the constituent maps
compatibly. Any morphism $f: A^\bl \ra B^\bl$ may be completed by a
mapping cone construction into a $3$-periodic sequence
\[ \ldots \ra A^\bl \ra B^\bl \ra C^\bl (f) \ra A^\bl[1] \ra \ldots\]
which becomes an exact sequence upon application of any functor
$\Hom_{D(\scrC)}(Z, -)$, and in particular upon taking cohomology.
This shows that although $D(\scrC)$ is not an abelian category, it is
what is known as a {\em triangulated} category.

\subsection{Derived functors}

We are particularly interested here in {\em morphisms} in the derived
category and in the way they compose.  They turn out to be
$\Ext$-groups, with composition being the Yoneda product. In other
words, the derived category is the place one should work if one wants
to view and compose elements of a cohomology group like morphisms ---
which is exactly what we propose to do to reformulate Kapranov's
construction of weight systems.

To explain this we need to consider derived functors. Suppose $F:
\scrC \ra \scrD$ is an additive functor between abelian
categories. Clearly it induces functors $\Ch(\scrC) \ra \Ch(\scrC)$
and $\Hty(\scrC) \ra \Hty(\scrD)$. But the obvious attempt to induce a
functor $D(F)\colon D(\scrC) \ra D(\scrD)$ between the derived
categories fails, because $F$ does {\em not} necessarily take
quasi-isomorphisms to quasi-isomorphisms. By considering mapping cones
one can see that this property is equivalent to $F$ taking all acyclic
complexes (those quasi-isomorphic to zero) to acyclic complexes, which
only holds for {\em exact} functors. To derive more general functors
we need to restrict the kinds of complex under consideration.

Recall that for any object $A \in \scrC$, the functor $\Hom_{\scrC}(-,
A): \scrC^{\op} \ra \Ab$ is {\em left-exact}, and that if it is also
{\em right-exact} then $A$ is called {\em injective}. Let
$\Inj(\scrC)$ denote the full subcategory of injective objects of
$\scrC$. If every object $A \in \scrC$ has an {\em injective
resolution} --- a quasi-isomorphism $A \ra I^\bl$ to a complex of
injective objects $I^\bl \in \Ch(\Inj(\scrC))$ --- then we say that
$\scrC$ {\em has enough injectives}.

Now any quasi-isomorphism out of an injective complex is a homotopy
equivalence; that is, we may construct its inverse in $\Hty(\scrC)$.
Consequently, any two injective resolutions of an object are
homotopy-equivalent, and if $\scrC$ has enough injectives then there
is an equivalence of categories between $\Hty(\Inj(\scrC))$ and
$D(\scrC)$. In this case one can define the {\em right-derived
functor} $RF$ of $F: \scrC \ra \scrD$ by just replacing $D(\scrC)$ by
$\Hty(\Inj(\scrC))$, applying $F$ to get to $\Hty(\scrD)$, and then
passing to the quotient $D(\scrD)$. Explicitly, if $B^\bl$ is an
object of $D(\scrC)$, one simply replaces it by an injective
resolution (well-defined up to homotopy equivalence) and applies $F$
to construct $RF(B^\bl)$.

The {\em classical derived functors} $R^iF$ associated to $F$ are just
the composites of the homology functors $h^i$ with $RF$. For the
functor $F=\Hom_{\scrC}(A, -)$ applied to an object $B \in \scrC$, we
have $R^iF(B)=\Ext_{\scrC}^i(A,B)$, because the procedure above agrees
with the traditional definition of the $\Ext$-groups: namely, take an
injective resolution of $B$, apply $\Hom_{\scrC}(A,-)$, and take
cohomology.

\subsection{Morphisms in the derived category}

Now we can explain the structure on the morphism sets in the derived
category which we need. The key fact is that for objects $A,B$ in
$\scrC$ we have
\[ \Hom_{D(\scrC)} (A, B[i]) =\Ext_{\scrC}^i(A,B).\]

Here is a sketch proof. First replace $B$ by an injective resolution
$I^\bl$, so that there is the isomorphism $\Hom_{D(\scrC)} (A, B[i])
=\Hom_{D(\scrC)} (A, I^\bl[i])$. Now elements of this latter group are
represented a priori by diagrams $A \ra C^\bl \leftarrow I^\bl$ whose
second map is a quasi-isomorphism; but because quasi-isomorphisms out
of an injective complex are invertible in $\Hty(\scrC)$, we only need
to look at actual homotopy classes of maps $A \ra I^\bl$. Finally,
chain homotopy classes of maps between chain complexes $A^\bl, B^\bl$
are given by the zeroth cohomology of the chain complex
$\Hom^\bl(A^\bl,B^\bl)$. Taking the shift into account gives the
result.

One important consequence is that we see that $\scrC$ is embedded in
$D(\scrC)$ as a full subcategory, because $\Hom_{D(\scrC)}(A,B) =
\Ext^0_{\scrC}(A,B) = \Hom_{\scrC}(A,B)$ for objects $A,B \in \scrC$.

There is a generalisation of the principle: for general objects of
$D(\scrC)$, complexes $A^\bl, B^\bl$, we have
\[ \Hom_{D(\scrC)} (A^\bl, B^\bl[i]) =\Ext_{\scrC}^i(A^\bl,B^\bl)\]
where the right-hand side is a ``hyperext'' group, computed by taking
an injective resolution of each of the terms of $B^*$, applying
$\Hom(A^\bl, -)$ and taking the total cohomology of the resulting
double complex.

It is also possible to show that the Yoneda product on $\Ext$ groups
of objects of $\scrC$
\[ \Ext_{\scrC}^i(A,B) \otimes \Ext_{\scrC}^j(B,C) \ra
\Ext_{\scrC}^{i+j}(A,C) \] corresponds to the composition of morphisms
$A \ra B[i]$ and $B \ra C[j]$ in $D(\scrC)$ that one gets after
applying the shift $[i]$ to the latter.


\subsection{Derived categories of coherent sheaves}
In this section we will look specifically at the case of the
derived category of $\OO_X$--modules on a complex manifold $X$.

Let $X$ be a finite-dimensional complex manifold $X$ and let
$\OO_X$ be the structure sheaf, that is its sheaf of germs of
holomorphic functions.  We are interested in the sheaves taking
into account the complex structure on $X$, these are the
\emph{sheaves of $\OO_X$-modules}, in other words sheaves $\scrE$
with a natural map of sheaves $\OO_X \otimes_\C \scrE \ra \scrE$.
The obvious example is the sheaf of germs of holomorphic sections
of a holomorphic vector bundle.  This example is {\em
locally-free} in the sense that any point has a neighbourhood $U$
over which the sections are isomorphic to the sheaf $\OO_U^{\oplus
k}$, for some $k$. The converse also holds: any locally-free sheaf
is the sheaf of sections of a holomorphic vector bundle.

We restrict the class of sheaves further by considering coherent
sheaves.  A {\em coherent sheaf} is a sheaf of $\OO_X$-modules
which is locally a quotient of a finite-rank locally-free sheaf.
On a smooth projective variety it actually has a global finite
resolution by locally-free sheaves. The coherent sheaves form an
abelian category and it is the \emph{bounded} derived category of
this that we refer to as the derived category of $X$ and denote
simply by $D(X)$.

We will use letters such as $E, F$ to denote locally-free sheaves,
script letters such as $\scrE, \scrF$ for general coherent sheaves,
and letters such as $A, B$ for typical objects in $D(X)$. The tangent
sheaf of $X$ will be written $T$ and its dual $T^*$ or $\Omega$. (We
abuse the star to indicate either the dual of a sheaf or a complex of
sheaves; it should be clear from the context which is intended.)

The {\em sheaf cohomology} groups $H^\gradestar (-)$ are the classical
derived functors of the global section functor $\Gamma = \Hom(\OO_X,
-)$, which takes sheaves to abelian groups. Thus, one computes
$H^\gradestar (\scrE)$ by taking an injective resolution of $\scrE$,
applying $\Gamma$ to obtain a chain complex of abelian groups, and
then taking the cohomology.  It can be helpful to have other points of
view: one can compute them using \v Cech cohomology, and for a
holomorphic vector bundle one can also think
differential-geometrically using the Dolbeault isomorphism $H^q(E)
\cong H_{\bar \partial}^{0,q}(E)$. For a compact complex manifold, all
cohomology groups are finite-dimensional.

In a similar vein, we can define the groups $\Ext^\gradestar (\scrE,
\scrF)$ by applying the classical derived functors of $\Hom(-,-)$ to
the pair $\scrE, \scrF$. They can also be computed by taking an
injective resolution of $\scrF$.

We can use the description of morphism sets from the previous section
to state the following result which is a key point for the
construction of the Rozansky-Witten weight systems: if $X$ is a
complex manifold and $\scrE$ is a coherent sheaf on $X$ then the sheaf
cohomology groups of $\scrE$ are expressible as morphism sets in the
derived category as follows:
\[H^q(\scrE) = \Ext^q(\OO_X, \scrE) = \Mor_{D(X)}(\OO_X,
\scrE[q]).\]

As another example of this logic, we can describe the cup product
\[ H^i(\textstyle{\bigwedge}^jT^\dualstar ) \otimes
H^k(\bigwedge^lT^\dualstar ) \ra 
H^{i+k}(\bigwedge^{j+l}T^\dualstar ) \] as the Yoneda product
operation which, given a pair of morphisms
\[ \OO_X \ra \textstyle\bigwedge^j{T}^\dualstar [i] \qquad \OO_X \ra
\bigwedge^l{ T}^\dualstar [k],\] applies the functor $-\otimes
\bigwedge^jT^\dualstar [i]$ to the second, composes the two, and
then performs exterior multiplication $\bigwedge^jT^\dualstar
\otimes \bigwedge^lT^\dualstar \ra \bigwedge^{j+l}T^\dualstar $.

Inside the category of coherent sheaves there is an internal
hom-functor: we can define $\sHom(\scrE,\scrF)$ to be the {\em sheaf}
of local homomorphisms $\scrE \ra \scrF$.  This has a right derived
functor which could be written $R\sHom(-,-)$ but which we will denote
for simplicity by $\sExt(-,-)$. The complex of sheaves
$\sExt({\mathcal E}, {\mathcal F})$ can be computed by taking a
locally-free resolution of $\scrE$ and applying $\sHom(-, \scrF)$.

The category of coherent sheaves is a tensor category under the
product $\otimes_{\OO_X}$, and the left derived functor of this
product equips the derived category $D(X)$ with the structure of a
symmetric tensor category. (We use an underline to distinguish the
derived functor $\underline{\otimes}$ from the underived $\otimes$
when applying it to complexes of sheaves, for which such a distinction
is necessary. But in a context where ``everything is derived'' we
often revert to the notation $\otimes$.) The identity object of $D(X)$
is the structure sheaf and the symmetry $\tau$ is the usual graded
symmetry for chain complexes. In fact there is a rigid structure: the
dual of an object $A$ is given by $A^*=\sExt(A, \OO_X)$. With this
definition, the double dual functor is canonically isomorphic to the
identity --- something which is not true for the naive (underived)
dualising functor $\Hom(-,\OO_X)$ defined on the category of coherent
sheaves.


\subsection{Standard operations with sheaves} 
\label{Section:SheafOperations}

For full details of all these operations and their relations, see
Kashiwara and Schapira \cite{KS}.

If $f\colon X \ra Y$ is a holomorphic map then there are induced
pushforward $f_*$ and pullback functors $f^*$ defined going between
the categories of coherent sheaves of $\OO_X$-modules and
$\OO_Y$-modules. These functor $f^*$ is left-adjoint to $f_*$, and
this relationship is preserved on the level of the derived category:
there are natural isomorpisms
\[   \Mor_{D(X)}(Lf^*A,B)\cong \Mor_{D(Y)}(A,Rf_*B).\]

One of the fundamental properties of the derived category of coherent
sheaves is that $Rf_*$ also has a right-adjoint $f^!\colon D(Y)\to
D(X)$, the Grothendieck-Verdier functor, so that
\[   \Mor_{D(X)}(B,f^!A)\cong \Mor_{D(Y)}(Rf_*B,A).\]
This functor $f^!$ can be defined as $Lf^*\Lotimes
Lf^*\omega_Y\Lotimes \omega_X [\dim X-\dim Y]$, where $\omega$
denotes the canonical line bundle, $\bigwedge^{\dim X}T^\dualstar$.

In fact these adjunctions hold ``internally'' in the derived category:
there are natural isomorphisms
\begin{align*}
   Rf_*\sExt(Lf^*A,B)&\cong \sExt(A,Rf_*B)\\
       \sExt(Rf^*B,A)&\cong Rf_*\sExt(B,f^!A).
\end{align*}
Other useful functorial identities are the tensoriality of the
pull-back
\[  Lf^*(A\Lotimes A')\cong Lf^*\Lotimes Lf^*A'\]
and the projection formula
\[   Rf_*(B\Lotimes Lf^*A)\cong Rf_*B \Lotimes A. \]


\subsection{Integral transforms}

Suppose we have two complex manifolds $X$ and $Y$. Then there is a
functor, {\em integral transform}, from $D(X\times Y)$ to the category
$\Fun(D(X),D(Y))$ of functors $D(X)\ra D(Y))$. Consider the diagram of
projections.
$$\begin{array}{ccc}
&X\times Y\\
&\llap{$\pi_X$}\swarrow\qquad\searrow\rlap{$\pi_Y$}\\
\qquad X&&Y\qquad.
\end{array}$$

We can view an object $P$ of $D(X \times Y)$ as a ``correspondence''
and define a functor $\widehat P\colon D(X) \ra D(Y)$, by pulling up
to $D(X\times Y)$, tensoring with $P$ and then pushing down to $Y$:
\[ \widehat P(A)=R\pi_{Y*}(\pi_X^*(A) \Lotimes P).\]
Here the pullback is exact and need not be derived. When $\widehat
P$ is an equivalence of categories, this is called a {\em
Fourier-Mukai transform}.

Moreover, a morphism $\Theta\colon P\to Q$ in $D(X\times Y)$ gives a
natural transformation $\widehat \Theta$ between the functors
$\widehat P$ and $\widehat Q$. Explicitly, we get for each object $A
\in D(X)$ a morphism
\[ \widehat\Theta_A:  R\pi_{Y*}(\pi_X^*(A) \Lotimes P) \ra
R\pi_{Y*}(\pi_X^*(A) \Lotimes Q) \] by applying the functors
$\pi_X^*(A) \Lotimes -$ and then $R\pi_{Y*}$ to the morphism
$\Theta$. So indeed we have a functor $\ \widehat \ : D(X \times
Y) \ra \Fun(D(X),D(Y))$. In what follows we will usually drop the hat
notation, using for example the same notation for morphisms in $D(X
\times Y)$ and their induced natural transformations.

Let us give some simple easily-checked examples of integral transforms
in the most important case, when $X=Y$. In this case we denote the two
projections by $\pi_1, \pi_2$, and we also consider the diagonal map
$\Delta: X \ra X \times X$.

The structure sheaf of the diagonal is an object $\Odel \in D(X)$,
given by the pushforward $\Odel = \Delta_* \OO_X$. This object gives
the identity functor $D(X) \ra D(X)$. If we look at the shifted
version $\Odel[n]$ it defines the shift functor $[n]: D(X) \ra D(X)$. 

We can define similarly define objects $T_\Delta = \Delta_* T$ and
$\Omega_\Delta = \Delta_* \Omega$ of $D(X \times X)$, which are
sheaves supported on the diagonal. It is easy to see that $\pi_1^* T
\otimes \Odel \cong T_\Delta$, and consequently (by means of the
projection formula) that $T_\Delta$ defines the ``tensor with $T$''
functor
\[ \id \otimes T : D(X) \ra D(X).\]

A little more subtly, for any object $A \in D(X)$ we can define the
{\em derived} pushforward $R\Delta_*A$ in $D(X\times X)$ and therefore
get an integral transform $D(X)\to D(X)$, which turns out to be just
the operation of derived tensor with $A$ (in $D(X)$). First notice
that there is an isomorphism of functors:
$R\Delta_*\cong\pi_1^*\Lotimes \Odel$.  This follows from a
straight-forward use of the projection formula:
\begin{align*}\pi^*_i(-)\Lotimes \Odel
  &\cong \pi^*_i(-)\Lotimes R\Delta_*\OO_X
  \cong R\Delta_*(\Delta^*\pi^*(-)\Lotimes \OO_X)\\
  &\cong
  R\Delta_*(\id_X(-)\Lotimes \OO_X)
  \cong R\Delta_*(-).
\end{align*}
Now we can see that indeed there is an isomorphism of functors
$\widehat{R\Delta_*A}\cong \id \Lotimes A$: just apply again the
projection formula:
\begin{align*}
  \widehat{R\Delta_*A}(-)
  &= R\pi_{2*}(\pi^*_1(-)\Lotimes R\Delta_*A)
  \cong R\pi_{2*}R\Delta_*(\Delta^*\pi^*_1(-)\Lotimes A)\\
  &\cong\id_X(\id_X(-)\Lotimes A)
  = (-)\Lotimes A.
\end{align*}

\pg


\section{The Atiyah class}


The construction we are interested in rests on the idea of the {\em
Atiyah class}, the sheaf-theoretic (and ultimately
derived-categorical) analogue of the curvature of a holomorphic
bundle. It is an extremely attractive and useful concept, so we devote
this section to a thorough explanation of its definition and
properties.

\subsection{The Atiyah class for vector bundles}

If $E$ is a holomorphic vector bundle on a complex manifold $X$ then
we can construct from a connection the curvature $1$-form $R_E$ used
in Section~\ref{Section:RWweightsystems}. Under the isomorphisms
\[  H^{0,1}_{\bar\partial}(E^\dualstar \otimes E \otimes T^\dualstar)
\cong H^1(E^\dualstar \otimes E \otimes T^\dualstar ) \cong \Ext^1(E
\otimes T, E) \] we can view it as a class $\alpha_E \in
\Ext^1(E\otimes T, E)$. Atiyah \cite{A} showed how to construct this
characteristic class in a purely sheaf-theoretic manner, giving it a
more canonical realisation.

One way to do this is as follows. If $E$ is a vector bundle, the
bundle of $1$-jets of $E$ is the sheaf $E \oplus E \otimes
\Omega$ with the twisted action of $\OO_X$ given by
\[
f. (s, t \otimes \theta) = (fs, ft \otimes \theta + s \otimes df)
\]
which describes first-order Taylor expansions of sections of $E$.
There is an exact sequence
\[
  0 \ra E \otimes \Omega \ra JE \ra E \ra 0,
\]
and the Atiyah class $\alpha_E \in \Ext^1(E, E \otimes \Omega )=
\Ext^1(E \otimes T, E)$ is defined to be the extension class.  The
extension class can be thought of as the obstruction to existence of a
section of the sequence -- in the case of a locally-free sheaf $E$,
such a thing would be a holomorphic connection on $E$ -- and may be
built in \v Cech cohomology using the differences between local
holomorphic splittings (which always exist). Another way to construct
it is to work purely homologically: tensoring the sequence with the
dual bundle $E^\dualstar $ gives another exact sequence whose
associated long exact sequence contains the map
\[
H^0(E^\dualstar  \otimes E)
 \stackrel{\delta}{\ra}
 H^1(E^\dualstar\otimes E \otimes \Omega ),
\]
and the Atiyah class is the image under $\delta$ of the identity
section of $\End(E)$.

The jet sequence/extension class definition also works for general
coherent sheaves $\scrE$, but the \v Cech cohomology representation is
more complicated in this case, since computing the relevant $\Ext$
group requires a resolution. We give a general recipe later in the
section.

\subsection{Properties of the Atiyah class}
The first important property of the Atiyah class we need is its {\em
naturality}. Suppose $f: E \ra F$ is a map of bundles, and regard each
of $\alpha_E$ and $\alpha_F$ as a morphism in $D(X)$. Then $f[1] \circ
\alpha_E = \alpha_F \circ (f \otimes \id_T)$, in other words the
diagram below commutes.
\[
 \begin{matrix} E \otimes T &\stackrel{\alpha_E}{\longrightarrow} &E[1]\\
 \llap{\scriptsize$f \otimes \id_T$} \downarrow &&
\llap{\scriptsize$f[1]$}\downarrow\\ F \otimes T
&\stackrel{\alpha_F}{\longrightarrow} &F[1]
\end{matrix}
\]

One way to prove this is to look at the long exact sequences
arising from $F^\dualstar $ tensor the jet sequence of $E$, from
$E^\dualstar $ tensor the jet sequence of $E$, and from
$F^\dualstar $ tensor the jet sequence of $F$. Since $f$ induces
maps of jet sequences, the latter two long exact sequences have
maps to the first one, and both identity elements map to $f\in
H^0(E^\dualstar \otimes F)$, proving a commutativity which when
written using morphisms is the one above.

A second important property is the behaviour under {\em tensor
product}.  One can show that the jet sequence for $E \otimes F$ is
the sum, in the sense of extensions, of the jet sequences $JE
\otimes F + E \otimes JF$. Thus the Atiyah class satisfies a
Leibniz rule as one might expect, which can be written sloppily as
\[ \alpha_{E \otimes F} = \alpha_E \otimes \id_F + \id_E \otimes
\alpha_F \] if we view this as an identity among morphisms $T
\otimes E \otimes F \ra E \otimes F[1]$.  Strictly speaking, some
permutations should be inserted to make this make sense, but there
are no sign problems until we deal with complexes of sheaves.

Finally, the Atiyah class of the tangent bundle $\alpha_T$ has a {\em
symmetry} and lies in fact in $\Ext^1(S^2 T, T)$. This corresponds in
differential geometry to the vanishing of the torsion (see
Section~\ref{Section:RWweightsystems}) and is explained elegantly by
Kapranov as follows. If $E$ is a sheaf on $X$ we can consider sheaves
of {\em $E$-torsors} over $X$, meaning sheaves whose local sections
are affine spaces modelled on the abelian group of local sections of
$E$. Thus, the sheaf $\Conn$ of local holomorphic connections on $X$
is a $T^\dualstar \otimes T^\dualstar \otimes T$-torsor. The torsion
defines a map from this sheaf to the sheaf of abelian groups
$\bigwedge^2T^\dualstar \otimes T$, and hence an exact sequence
\[ 0 \ra \Conn_{tf} \ra \Conn \ra \textstyle\bigwedge^2T^\dualstar
\otimes T \ra 0,\] where the first term is the sheaf of torsion-free
connections, a torsor over $S^2T^\dualstar \otimes T$.  Each torsor
defines an obstruction element in $H^1$ of its appropriate model
sheaf. These elements are related by the long exact sequence arising
from
\[
0 \ra S^2T^\dualstar  \otimes T
  \ra T^\dualstar  \otimes T^\dualstar  \otimes T
  \ra \textstyle\bigwedge^2T^\dualstar  \otimes T \ra 0,
\]
so the fact that $\bigwedge^2T^\dualstar  \otimes T$ is a trivial
torsor means it represents the trivial element, and so the Atiyah
class comes from a symmetric element.


\subsection{Functorial definition of the Atiyah class}

We will need to extend the definition of the Atiyah class from bundles
to general objects of the derived category. For each object $A$, we
would like an element $\alpha_A \in \Ext^1(A \otimes T, A)$, or
equivalently a morphism in $D(X)$
\[ \alpha_A : A \otimes T \ra A[1]. \]
(Remark: here and subsequently, if a tensor product here is {\em
obviously} derived, as for example when we are dealing with $D(X)$ as
a tensor category, we do not distinguish it by an underline.) The
naturality square from the previous section suggests that such
morphisms should form the components of a natural transformation
\[ \alpha : \id \otimes T \ra \id[1]\]
and this is exactly what we establish below. One way to do this is to
build explicitly a representative for complexes of locally-free
sheaves, starting from the above version for single sheaves. (We will
shortly give a \v Cech description of the Atiyah class which could be
used to do this.) But there is a far more elegant way to construct
$\alpha$ directly.

Consider the product $X \times X$, with the two projections $\pi_1,
\pi_2$ and the diagonal $\Delta \subseteq X \times X$. Associated to
$\Delta$ is the ideal sheaf $\scrI_\Delta$ of holomorphic functions on
$X\times X$ vanishing on $\Delta$, and there is an exact sequence
\[
0 \ra \scrI_\Delta/\scrI_\Delta^2
  \ra \OO_{X \times X}/\scrI_\Delta^2
  \ra \OO_{X \times X}/\scrI_\Delta
  \ra 0,
\]
whose three terms are identifiable respectively as: the cotangent
sheaf $\Omega_\Delta \cong \pi_1^*\Omega \otimes \Odel$ of $\Delta$;
the structure sheaf of the first infinitesimal neighbourhood of
$\Delta$; and the structure sheaf $\Odel$ of $\Delta$. This sequence
defines an extension class
\[ \alpha \in \Ext^1_{X \times X}(\Odel, \Omega_\Delta) = \Hom_{D(X
  \times X)}(\Odel, \Omega_\Delta[1])  =
\Hom_{D(X \times X)}(T_\Delta, \Odel[1]).\] Therefore it gives, by
integral transform, a natural transformation
\[ \alpha: \id \otimes T \ra \id[1]\]
between the ``tensor with $T$'' and ``shift by $1$'' functors, as
required.

We think of the morphism $\alpha \in \Hom_{D(X\times
X)}(T_\Delta,\Odel[1])$ as the ``universal Atiyah class'' for $X$. We
should check that from it we can indeed recapture the earlier
definition of the Atiyah class, in the case when $A=E$ is a single
locally-free sheaf. For this we only need to observe that if we apply
the functor $R\pi_{2*}(\pi_1^*E\Lotimes {-})$ to the sequence
\[ 0 \ra \Omega_\Delta \ra \OO_{X \times X}/\scrI_\Delta^2 \ra \Odel
\ra 0 \] we get the jet sequence for $E$. Consequently, the universal
Atiyah class $\alpha: \Odel \ra \Omega_\Delta[1]$, which extends the
first sequence into a distinguished triangle, is sent to the Atiyah
class $\alpha_E$, which extends the latter to a distinguished
triangle. (Recall that derived functors preserve distinguished
triangles.)

The properties of the Atiyah class that we observed for bundles still
hold in this more general context. The naturality follows
automatically from the construction via the universal class $\alpha$.
The Leibniz tensor product rule still holds for this generalised
Atiyah class, with the symmetry $\tau$ taking care of the signs, and
the symmetry property, which is special to the tangent sheaf $T$, is
unchanged.

\subsection{Explicit representation of the Atiyah class.}

Although we have tried to define the Atiyah class in the most elegant
way possible, the abstract definition sometimes needs to be
supplemented by a way of actually calculating it in examples. We give
here a rather long exposition of how to do this, and most readers
should probably ignore it, since in fact we only need this result at
one point in section 7. 

Recall the construction of the connecting homomorphism in the long
exact sequence of cohomology of sheaves. If 
\[ 0 \ra E \ra F \ra G \ra 0\]
is an exact sequence then we take injective resolutions of these three
sheaves, obtaining an exact sequence of complexes of sheaves
\[ 0 \ra I^* \ra J^* \ra K^* \ra 0,\]
apply the section functor $\Gamma$ to get an exact sequence of
complexes of abelian groups
\[ 0 \ra \Gamma(I^*) \ra \Gamma(J^*) \ra \Gamma(K^*) \ra 0.\]
Then the standard Snake lemma construction defines the coboundaries
\[ H^i(G) \ra H^{i+1}(E).\]

A much more tangible version is obtained by using \v Cech complexes
instead. Fix some good cover of $X$ and let $C^*(E) = C^0(E) \ra
C^1(E) \ra \cdots$ be the associated \v Cech
complex. Then we have an exact sequence of complexes of abelian groups
\[ 0 \ra C^*(E) \ra C^*(F) \ra C^*(G) \ra 0.\]
Via the usual double complex proof (look at $C^*(I^*)$ where $I^*$ is
an injective resolution of $E$) we know that \v Cech and sheaf cohomology
are isomorphic, and since this isomorphism is functorial, it follows
that the connecting homomorphisms coming from this exact sequence
coincide with the ones from the first construction via injective
resolutions. 

If we now generalise to the case of hypercohomology, where we are
computing the cohomology of a complex $E^*$ of sheaves, the same
arguments go through ``with an additional index'' as follows. By
definition we compute the hypercohomology by forming a double complex
of injective resolutions of the sheaves of $E^*$, and taking its total
cohomology. But an alternative method of computation is to form a
double complex $C^*(E^*)$ of \v Cech complexes of the sheaves making up
$E^*$, and to take its total cohomology. (The proof that these two
methods are functorially isomorphic uses a triple complex!) The
connecting homomorphisms in the long exact sequence of hypercohomology
come from the Snake lemma applied to a three-term sequence of double
complexes (by using their total differentials), and we can construct them 
similarly in \v Cech cohomology.

These principles give us a way to write down representatives of the
Atiyah class. We deal first with the case of a single locally-free
sheaf $E$ (this is very easy but it is a great help in explaining
the more complicated case) and then with a complex of locally-free
sheaves. As a general object in the derived category $D(X)$ is
quasi-isomorphic to such a complex, this is all we ever need.

Recall that $\alpha_E \in \Ext^1(E, E \otimes \Omega)$ is
an obstruction class: it is the image of the identity under the
connecting homomorphism
\[ \Ext^0(E, E) \ra \Ext^1(E, E \otimes \Omega)\] coming from the
long exact sequence of classical derived funtcors $\Ext^i(E, -)$
applied to the jet exact sequence
\[ 0 \ra E \otimes \Omega \ra JE \ra E \ra 0.\]

Let's view $\Ext^i$ in this context as the composite of functors $h^i
\circ R\Hom = h^i \circ R\Gamma \circ R\sHom$ ($h^i$ denotes, as
usual, the $i$th cohomology of a complex). Since $E$ is locally-free,
the functor $\sHom(E, -)$ is exact, and hence $R\sHom(E, E \otimes
\Omega)$ is simply the sheaf $\sHom(E, E \otimes \Omega)$. Therefore
\[ \Ext^i(E, E \otimes \Omega) = H^i(\sHom(E, E \otimes \Omega))\]
can be viewed as simply a sheaf cohomology group.

The relevant connecting homomorphism $\delta$ comes from the diagram
of \v Cech cochain groups
\[\begin{matrix}
0 &\ra& C^0(\sHom(E, E \otimes \Omega)) &\ra& C^0(\sHom(E, JE)) &\ra&
C^0(\sHom(E,E)) &\ra& 0\\
&&\downarrow&&\downarrow&&\downarrow&&\\
0 &\ra& C^1(\sHom(E, E \otimes \Omega)) &\ra& C^1(\sHom(E, JE)) &\ra&
C^1(\sHom(E,E)) &\ra& 0.\end{matrix}\]
Let us compute $\delta(1)$. Fix a good cover
$\{U_i\}$ of $X$ over which $E$ is locally trivial. Begin with the
0-cochain $\{1_i\}_i$ in the top right ($1_i$ represents the identity
of $E\vert_{U_i}$). We lift this to a cochain in
the top middle. Since $JE = E \oplus E \otimes \Omega$ and the top
right map is just projection to $E$, the lift must be of the form
$\{1_i \oplus \nabla_i\}_i$, where $\nabla_i: E\vert_{U_i} \ra (E
\otimes \Omega)\vert_{U_i}$ satisfies (for $f \in \Gamma(\OO_X, U_i)$
and $s \in \Gamma(E, U_i)$)
\[ \nabla_i(f.s)= f.\nabla_i s + f. ds\]
and is therefore a {\em connection} on $E$ over $U_i$. Since $E$ is
trivial on each $U_i$, such a thing exists.

Now applying the \v Cech coboundary and lifting to the bottom left
corner, we end up with the cochain $\{\nabla_i -\nabla_j\}_{ij} \in
C^1(\sHom(E, E \otimes \Omega))$. Clearly we have recovered the fact
that the Atiyah class is the obstruction to existence of a global
holomorphic connection.

When $(E^*, \partial)$ is a {\em complex} of locally-free sheaves we
modify this construction as follows. Once more we view $\Ext^i$ as the
composite of functors $h^i \circ R\Hom = h^i \circ R\Gamma \circ
R\sHom$. Since $E^* \otimes \Omega$ is locally-free, the functor
$\sHom(E^* \otimes \Omega, -)$ (taking complexes of sheaves to
complexes of sheaves) is exact, and hence $R\sHom(E^*, E^* \otimes
\Omega)$ is simply the complex of sheaves $\sHom(E^*, E^* \otimes
\Omega)$. Therefore
\[ \Ext^i(E^*, E^* \otimes \Omega) = H^i(\sHom(E^*, E^* \otimes \Omega)).\]
is just a hypercohomology group. It can be computed from the total
cohomology of the double complex $C^*(\sHom(E^*, E^* \otimes \Omega))$
(the \v Cech complex of a complex of locally-free sheaves).

To compute the connecting homomorphism we use the analogue of the
diagram above. This time the groups are the total cochain groups of
double complexes and the vertical coboundary maps are the total
differentials in these double complexes, namely $d+(-1)^p\partial$,
where $d$ is the \v Cech differential and $\partial$ the differential
on the complex $E^*$. We begin with the collection of identity maps
$\{1^i_j\}\in C^0(\sHom^0(E^*, E^*))$. (Here the lower index denotes
the set of the cover and the upper one the position in the complex, so
that $1^i_j$ is the identity $E^i(U_j) \ra E^i(U_j)$.) For each sheaf
$E^i$ and open set of the cover $U_j$ we pick a local connection
$\nabla^i_j$ so that the lift of the identity is $\{1^i_j \oplus
\nabla^i_j\} \in C^0(\sHom^0(E^*, JE^*))$. Now apply the vertical
coboundary and lift into the bottom left corner, namely
\[ C^0(\sHom^1(E^*, E^* \otimes \Omega)) \oplus C^1(\sHom^0(E^*, E^*
\otimes \Omega)). \]
The result is that the Atiyah class is represented by 
\[ \{\partial \nabla^i_j -
\nabla^i_j \partial\}_j \oplus \{ \nabla^i_j - \nabla^i_k\}_{jk} \in
C^0(\sHom^1(E^*, E^* \otimes \Omega)) \oplus C^1(\sHom^0(E^*, E^*
\otimes \Omega)). \]

Of course in the special extremal case that $E^*$ is a single sheaf
the first term drops out and we get back the representative we already
computed. In the other extremal case where the $E^*$ are globally
trivial (for example on an affine space), the second term drops out
and we just have $\partial \nabla - \nabla \partial$ as
representative. We refer to this statement as {\em Markarian's lemma
$1$}, since it comes from his paper \cite{Ma} (in which it is an
exercise for the reader). 

Finally we observe that to compute the Atiyah class for an arbitrary
(not locally-free) sheaf or complex of sheaves $\scrE^*$ we can just
resolve first by a (double) complex of locally-free ones and then use
the above method to obtain a representative of the Atiyah class.

\subsection{Final comments on the Atiyah class}

There are a few further comments we will need soon.

{\em The Atiyah class of the diagonal}.
Recall that the universal Atiyah class
\[ \alpha \in \Ext^1_{X \times X}(T_\Delta, \Odel)\]
comes (after taking an adjoint) from the infinitesimal neighbourhood
sequence
\[ 0 \ra \Omega_\Delta \ra \OO_{X \times
X}/\scrI_\Delta^2 \ra \Odel \ra 0.\]

This morphism is very closely related to the Atiyah class of $\Odel$
itself, which lies in $\Ext^1_{X \times X}(\Odel \otimes T_{X\times
X}, \Odel)$ and is (by definition) the extension class of the jet
sequence
\[ 0 \ra \Odel \otimes \Omega_{X \times X} \ra J(\Odel) \ra
\Odel \ra 0.\]
We can decompose $\Omega_{X \times X} = \pi_1^*\Omega_X \oplus
\pi_2^* \Omega_X$ and therefore identify 
\begin{eqnarray*} \Ext^1_{X \times X}(\Odel \otimes
T_{X\times X}, \Odel) &=& \Ext^1_{X \times X}(\Odel \otimes \pi_1^*T,
\Odel) \oplus \Ext^1_{X \times X}(\Odel \otimes \pi_2^*T, \Odel)\\ &=&
\Ext^1_{X \times X}(T_\Delta, \Odel) \oplus \Ext^1_{X \times
X}(T_\Delta, \Odel).\end{eqnarray*} It is easy to check that the jet
sequence is a Baer sum of two copies of the infinitesimal
neighbourhood sequence and hence that under this identification the
Atiyah class $\alpha_{\Odel}$ is equal to the sum $\alpha \oplus
\alpha$ of two copies of the universal Atiyah class. This remark will
be important in understanding the STU relation for the universal
Atiyah class $\alpha$. 

{\em Locality}.  From abstract functoriality, or directly from the
local representation of the Atiyah class, one can see the following
locality property: if $E$ is an object of $D(X)$ and $U \subseteq X$
is an open set, then the diagram
\[\begin{matrix} E &\stackrel{\alpha_X}{\ra}& E \otimes \Omega_X \\
\downarrow && \downarrow\\ i^*E &\stackrel{\alpha_U}{\ra}& i^*E \otimes
\Omega_U \end{matrix}\] commutes. 

{\em Functoriality under pullback}. 
A final property we need is about pullbacks of the Atiyah class. This
is that the diagram
\[\begin{matrix} f^*E \otimes T_Y &\stackrel{\alpha_{f^*E}}{\ra}& f^*E[1]\\
\downarrow&&\downarrow \\f^*E \otimes f^*T_X
&\stackrel{f^*(\alpha_E)}{\ra}& f^*E[1],\end{matrix}\] where the
left-hand downward map is $\id \otimes df$, commutes. As an example,
consider the Atiyah class $\alpha_{T_{X \times X}} \in \Ext^1(T_{X
\times X} \otimes T_{X \times X}, T_{X \times X})$. It can obviously
be decomposed into two pieces via the usual splitting, with the first
living in $\Ext^1(\pi_1^* T \otimes \pi_1^*T, \pi_1^*T)$. This piece
equals $\pi_1^*\alpha_T$, by the above naturality.

\pg


\section{Rozansky-Witten weight systems revisited}


In this section we bring together the abstract nonsense of the
previous three sections and show how it provides an elegant
formulation of Rozansky-Witten weight systems. 


\subsection{The Lie algebra object of a complex manifold}
The first main theorem of the paper is the following interpretation of
a complex manifold as ``being'' in some sense a Lie algebra. 

\begin{thm}
Suppose $X$ is a complex manifold. Then the shifted tangent sheaf
$T[-1]$ is a Lie algebra object in the derived category $D(X)$;
furthermore, every object in $D(X)$ is canonically a module over
$T[-1]$, and every morphism in $D(X)$ is a module map.
\end{thm}
\begin{proof}
We just need to define the structure morphisms and check the
identities for them. To obtain the bracket, start with the Atiyah
class of $T$, viewed as a morphism $\alpha_T:  T \otimes T \ra
T[1]$. Now apply an additional shift by $[-2]$ to each side and the
result is the bracket
\[ T[-1] \otimes T[-1] \ra T[-1]. \]
The module action for any object $A \in D(X)$, likewise, is just
obtained by shifting the Atiyah class $\alpha_A$ by $[-1]$; it is a
morphism
\[ A \otimes T[-1] \ra A. \]
Skew-symmetry of the bracket comes because the {\em unshifted} Atiyah
class is symmetric, and the shifts of $[-1]$ switch the parity. The
Jacobi (IHX) identity and module (STU) identity are just the fact that
the two morphisms above are {\em invariant} under the action of the
Atiyah class, which is a consequence of its naturality. Explicitly,
for the STU case: consider the morphism $\alpha_A:  A \otimes T \ra
A[1]$. This commutes with taking Atiyah classes on each
side, according to the diagram
\[
\begin{matrix}
 (A \otimes T)  \otimes T
     &\stackrel{\alpha_{A\otimes T}}{\longrightarrow} &(A \otimes T)[1]\\
  \llap{\scriptsize$\alpha_A \otimes \id$} \downarrow &
       & \llap{\scriptsize$\alpha_A[1]$}\downarrow\\
  A \otimes T[1] &\stackrel{\alpha_A[1]}{\longrightarrow} &A[2]
\end{matrix}.
\]
Using the Leibniz rule to evaluate the top line, and putting in
the shifts (this affects the signs a little) gives the STU
relation. It makes more sense with pictures; naturality and the
Leibniz rule for the the Atiyah class amount to the identity
\[ \vpic{xibottom} = \vpic{xitop}\]
for {\em any} (boxed) morphism between two tensor products of objects in
$D(X)$. Applying this naturality to $\alpha_A[-1]$ (the case where the
box is actually a trivalent vertex) gives the familiar
\[ \vpic{n1} + \vpic{n2} = \vpic{n3},\]
and the IHX relation is the special case $A=T[-1]$. Note once more the
way that the signs are {\em locally} compatible with the ``correct''
vertex-oriented orientations of graphs.
\end{proof}
There are two things to note here.  Firstly, the Lie algebra
structure on $X$ is nilpotent: an $n$-fold composition defined
using the Atiyah class (corresponding pictorially to $n$-vertex
tree with $n+1$ inputs and $1$ output) lies in $\Ext^n(T^{\otimes
n+1}, T)$, which is zero for $n$ larger than the complex dimension
of $X$.

Secondly, the bracket $\alpha: T \otimes T \ra T[1]$ induces by
composition a bracket on the vector space $\bigoplus_n\Mor(\OO_X,
T[n]) = H^\gradestar (T)$. One might regard this as the simplest
``computable''manifestation of the Atiyah class, but unfortunately it
is zero. This is because the composite of any two elements may be
pictured as below, and sliding the trivalent vertex down past one of
the boxes creates the Atiyah class of the structure sheaf (the unit
object) which is zero.
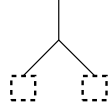
\[ \vpic{blobs} \]

\subsection{Metric Lie algebras from complex symplectic manifolds}

The second main theorem is the similar interpretation of complex {\em
symplectic} manifolds as ``being'' {\em metric} Lie algebras.

A complex symplectic form $\omega \in H^0(T^\dualstar \otimes
T^\dualstar )$ may be rewritten as a morphism $T \otimes T \ra \OO_X$
and then, by shifting by $[-2]$, as a {\em symmetric} morphism
\[ T[-1] \otimes T[-1] \ra \OO_X[-2].\]

As a morphism in $D(X)$, this is invariant under the action of $T[-1]$
and so satisfies the identities stated in section 3. However, it
cannot quite be regarded as metric on $T[-1]$ because of the shifts
$[-2]$ appearing on the right-hand side. A metric on $L$ is meant to
be a morphism $L \otimes L \ra 1$, which in our case would be a
morphism $T[-1] \otimes T[-1] \ra \OO_X$, without the shift.  To
handle this difficulty we alter $D(X)$ into a category $\tilde D(X)$:
we define it to have the same objects as $D(X)$ but redefine the space
of morphisms $A \ra B$ to be the graded vector space $\Ext^\gradestar
(A,B)$ instead of just $\Ext^0(A,B)$. Composition of morphisms is
defined in the obvious way and is graded bilinear. After this
extension, the above shifts cease to cause problems. In summary:

\begin{thm}
If $X$ is a complex symplectic manifold then $T[-1]$ is a metric Lie
algebra in the extended derived category $\tilde D(X)$, and $\tilde
D(X)$ is a module category over $T[-1]$.
\end{thm}

Consequently we can apply the general categorical construction of
weight systems from section 3:

\begin{thm}
If $X$ is a complex symplectic manifold and $A$ is an object of
$D(X)$ then there is a weight system
\[ w_{X,A}: \A \ra H^\gradestar(\OO_X).\]
\end{thm}

\begin{remk}
A different way to define $\tilde D(X)$ is as follows. Embed $D(X)$ in
the derived category $D^u(X)$ of unbounded complexes using the functor
$i=\bigoplus_{n\in\N} [n]$: an object $A \in D(X)$ is sent to $i(A) =
\bigoplus_{n\in\N} A[n]$. The set of morphisms $i(A) \ra i(B)$ in
$D^u(X)$ is the rather large $\bigoplus_{m,n\in \N} \Mor_{D(X)}(A[m],
B[n])$, but the shift functor $[1]$ acts on this space, and the set of
morphisms which {\em commute} with this action is the very reasonable
\[\bigoplus_{n\in \N} \Mor_{D(X)}(A, B[n]) =\Ext^\gradestar (A, B).\]
By this procedure of essentially looking at the $[1]$-invariant
subcategory of $D^u(X)$, we define $\tilde D(X)$.

This odd-looking construction has the advantage of being exactly
parallel to the procedure of replacing the category of
finite-dimensional complex $\g$-modules with modules over
$\C[[\hbar]]$. One replaces every space $V$ with the graded space
$\bigoplus_{n\in\N} V.\hbar^n (=V \otimes_{\C} \C[[\hbar]])$, and uses
only the maps of $\C[[\hbar]]$ modules, that is the
$\hbar$-equivariant linear maps between these. 

This construction is necessary in the theory of Vassiliev invariants
if we want to obtain from a Lie algebra $\g$ a weight system defined
on the graded completion of $\A$. In order to avoid convergence
problems we have to multiply the Casimir element by an indeterminate
$\hbar^2$ and the metric by $\hbar^{-2}$ to obtain weight systems
$\hat \A \ra \Q[[\hbar^2]]$. Salvaging some of the grading in this way
is absolutely essential to the correspondence between the Kontsevich
integral and invariants coming from quantum groups, and to the
deformation of the category of representations of $\g$ via the
Knizhnik-Zamolodchikov equation. In section 9 we will see the parallel
deformation for $\tilde D(X)$.
\end{remk}

\begin{remk}
There are weaker geometrical structures we could consider. If $X$ is a
{\em holomorphic Casimir manifold}, possessing a holomorphic bivector
$w\in H^0(\Lambda^2T)$ (not required to be non-degenerate) then
$T[-1]$ is a Casimir Lie algebra in $\tilde D(X)$, in an analogous
way. The Casimir is the symmetric morphism $\OO_X \ra (T[-1] \otimes
T[-1])[2]$. We can also formulate the even weaker analogue of a vector
space with a classical $r$-matrix too: this is a complex manifold $X$
with (for example) a sheaf $\scrE$ and an element $r \in
\Ext^\gradestar (\scrE \Lotimes \scrE, \scrE \Lotimes \scrE)$
satisfying the 4T relation of Vassiliev theory. But this is probably
not very useful.
\end{remk}

\subsection{The STU relation for the universal Atiyah class}

We've seen that for any object $A \in D(X)$, its Atiyah class
$\alpha_A: A \otimes T[-1] \ra A$ together with that of the tangent
sheaf $\alpha_T: T[-1] \otimes T[-1] \ra T[-1]$ satisfy the STU
relation, which can be written non-pictorially as 
\[ \alpha_A \circ \alpha_T = [\alpha_A, \alpha_A] \in \Hom_{D(X)}(A \otimes
T[-1] \otimes T[-1], A).\] This strongly suggests that the universal
Atiyah class morphism $\alpha: T_\Delta[-1] \ra \OO_\Delta$ in
$D(X\times X)$, together with the pullback $\pi^*(\alpha_T):
\pi^*T[-1] \otimes \pi^*T[-1] \ra \pi^*T[-1]$ (the sources of the above
morphisms), should satisfy the corresponding ``universal'' relation
\[ \alpha \circ \pi^*(\alpha_T) = [\alpha, \alpha] \in \Hom_{D(X
  \times X)}(\OO_\Delta \otimes \pi^*T[-1] \otimes \pi^*T[-1],
  \OO_\Delta).\] This is in fact the case because of the relation
  between $\alpha$ and the Atiyah class of $\Odel$. Certainly we have
  the identity
\[ [\alpha_{\Odel}, \alpha_{\Odel}] = \alpha_{\Odel} \circ \alpha_{T_{X
\times X}},\] and if we extract the first component parts of these
identities under the usual splitting of $T_{X \times X}$ we get the
desired equality.


\pg

\section{The symmetric and universal enveloping algebras of $T[-1]$}

Let $\g$ be a Lie algebra. The symmetric algebra $S(\g)$ and the
universal enveloping algebra $U(\g)$ are defined as quotients of the
tensor algebra $\TT (\g)=\bigoplus\g^{\otimes n}$:
\[
S(\g):=\TT(\g)/\langle x\otimes y-y\otimes x\rangle;\quad
U(\g):=\TT(\g)/\langle[x,y]-x\otimes y+y\otimes x\rangle.
\]
Each inherits an associative algebra structure and $\g$-module
structure from the tensor algebra; the symmetric algebra also inherits
a grading. The universal enveloping algebra has a universal property
for Lie algebra homomorphisms from $\g$ into associative algebras, and
the representation theory of $U(\g)$ coincides with that of $\g$.

Via symmetrization there is a splitting $S(\g)\hookrightarrow\TT(\g)$
and composing this with the quotient map $\TT(\g)\twoheadrightarrow
U(\g)$ gives a vector space isomorphism called the
Poincar\'e-Birkhoff-Witt map:
\[
 \PBW\colon S(\g)\to U(\g).
\]
This is a $\g$-module map, so it induces a vector space isomorphism
on the invariant parts:
\[
 \PBW\colon S(\g)^\g \cong U(\g)^\g.
\]

We have seen that the object $L=T[-1]$ is a Lie algebra for any
complex manifold and that $D(X)$ is a category of modules over $L$. We
now pursue this analogy further: we construct objects $S$ and $U$, the
symmetric and universal enveloping algebras of $L$, and a PBW
isomorphism between them. The third main theorem of the paper is:

\begin{thm}
The object $S = \bigoplus (\bigwedge^kT)[-k]$ is the symmetric algebra of
$L$, while $U = \pi_*\sExt(\Odel, \Odel)$ is its universal enveloping
algebra. These objects satsify the expected universal properties, and
$U$ acts on all objects of $D(X)$ compatibly with the action of
$L$. There is a morphism $\PBW : S \ra U$, the PBW morphism,
which is an isomorphism of objects (but not of algebras) in $D(X)$.
\end{thm}

The construction of $S$ is straightforward, but verifying the
properties of $U$ is quite difficult, and relies on some ideas of
Nikita Markarian \cite{Ma}. Andrei C\u ald\u araru \cite{CaldararuII}
independently explored similar ideas, to a different purpose, and very
recently Ajay Ramadoss \cite{Ram} studies a similar problem.

To simplify notation, in this section {\em all functors will be
derived}, so $\otimes$ means $\Lotimes$, $f^*$ means $Lf^*$, and $f_*$
means $Rf_*$. (With this convention we could write $\sHom$ for
$\sExt$, but we won't.) We will also write simply $\pi$ for the
projection $\pi_1: X \times X \ra X$.

\subsection{The symmetric algebra}

The symmetric power $S^k(T[-1])$ is actually the object
$(\bigwedge^kT)[-k]$, because the shift $[-1]$ changes the parity of
the flip map $\tau$ in $D(X)$ and therefore changes symmetrisation to
antisymmetrisation. Thus, the symmetric algebra of $T[-1]$ is the
object
\[ S = \bigoplus (\textstyle\bigwedge^kT)[-k].\]
It is a finite sum and is equipped with the commutative algebra
structure induced by exterior multiplication. It is easy to see that
it is category-theoretically the symmetric algebra $S(L)$ of
$L$. Firstly there is a canonical map $L \ra S$. Secondly, given any
map from $L$ to a commutative algebra object $A$, we get a lift $S \ra
A$ by symmetrisation (view $\bigwedge^kT$ as a subsheaf of the tensor
sheaf) followed by multiplication in the normal way. This gives a
commutative algebra homomorphism, uniquely determined by the original
$L \ra A$.

\subsection{The universal enveloping algebra: plan of attack}

The usual construction in the category of vector spaces builds $U(\g)$
as a quotient of $T(\g)$. We cannot do the same construction in $D(X)$
because it is not an abelian category, merely triangulated. In any
case, we want to have a reasonable description of the object $U$, not
simply an abstract definition as a quotient. Our definition $U =
\pi_*\sExt(\Odel, \Odel)$ is quite explicit, but it is unfortunately
relatively hard to show that it really is the universal enveloping
algebra of $L$, in the sense of category theory. (While this isn't
really essential to our study of Rozansky-Witten invariants, it is
worth establishing in its own right and is conceptually important in
studying the TQFT. )

Here are the steps we must take to prove the theorem.

1. Show that $U$ is an associative algebra object 

2. Construct a natural map $L \ra U$ which is a Lie algebra
   homomorphism (with respect to the commutator bracket on $U$).

3. Show that the universal property holds: every Lie algebra morphism
   $L \ra A$ for some other associative algebra $A$ extends (under
   $L\ra U$) to an associative algebra morphism $U \ra A$.

4. Construct a map $S \ra U$ which is an isomorphism of objects in
   $D(X)$.

5. Show that $U$ acts on all objects in $D(X)$, compatibly (under
   $L\ra U$) with the action of $L$.

It is relatively straightforward to perform steps 1, 2 and 5 and this
is handled in the next subsection.

Step 4, the construction of the PBW morphism, was done by Markarian
\cite{Ma} in lemma 1 (``proof: left to reader.'') and
definition-proposition 1 (``proof: local check is enough.''). Not
being experts, we didn't find these exercises at all trivial, so we
worked out the details, the first in section 5 (the local
representation of the Atiyah class) and the second below. (Although
these are a bit long-winded, we felt it would be useful to supply
details as an aid to anyone else who has tried to understand
Markarian's paper.)

Step 3 is the most frustrating step: we know of no direct way of
constructing the requisite maps $U \ra A$. So instead we fall back on
a rather abstract method of proof which relies on steps 1,2 and 4 and
a theorem of Hinich and Vaintrob. This is contained in the penultimate
subsection, after which there are some further remarks on the
structure of $S$ and $U$.

\subsection{The construction of $U$.}

We define $U = \pi_* \sExt(\Odel, \Odel)$, an object of $D(X)$. 

\noindent{\em Step 1}.  This object $U$ is an associative algebra in
$D(X)$. To see this, first observe that $A=\sExt(\Odel, \Odel)$ is an
associative algebra in $D(X \times X)$. If we apply the pushforward to
the multiplication map $A \otimes A \ra A$ then we get a map
\[ \pi_*(A \otimes A) \ra \pi_*A\]
which is not quite what we want. However, there is a natural map (the
adjunction unit)
\[ \pi^*\pi_*A \ra A\]
in $D(X \times X)$ corresponding to the identity under the adjunction
isomorphism
\[ \Hom_{D(X\times X)}(\pi^*\pi_*A, A) \cong  \Hom_{D(X)}(\pi_*A, \pi_*A)\]
and if we tensor this with itself we get a map
\[  \pi^*( \pi_*A \otimes \pi_*A) = \pi^*\pi_*A \otimes
\pi^*\pi_*A\ra A \otimes A\] whereupon the adjunction isomorphism
(reversed) gives us a map
\[ \pi_*A \otimes \pi_*A \ra \pi_*(A \otimes A).\]
Precomposing with this gives us the required multiplication $U \otimes
U \ra U$. It is straightforward to check that it is still associative
and unital.

\noindent{\em Step 2}.
Next, we define the canonical Lie algebra homomorphism $\gamma:L \ra
U$. Consider the universal Atiyah class morphism
\[ \alpha \in \Hom_{D(X\times X)}(T_\Delta, \Odel[1]) \cong
\Hom_{D(X\times X)}(\pi^*T[-1] \otimes \Odel, \Odel)\]
in the adjoint form (moving $\Odel$ to the RHS) 
\[ \alpha \in \Hom_{D(X\times X)}(\pi_1^* T[-1], \sExt(\Odel, \Odel) )\]
and apply the adjunction
\[ \Hom_{D(X\times X)}(\pi_1^* T[-1], \sExt(\Odel, \Odel) ) \cong \Hom_{D(X)}(T[-1],
\pi_* \sExt(\Odel, \Odel))\] to get the required map $\gamma: T[-1]
\ra \pi_* \sExt(\Odel, \Odel))$.

We must show that this is a morphism of Lie algebras when $U$ is given
the commutator bracket, that is that the diagram
\[ \begin{matrix} T[-1] \otimes T[-1] &\stackrel{\gamma
    \otimes \gamma}{\ra} & U \otimes U\\
\alpha_T \downarrow &&\downarrow [,]\\
T[-1] &\stackrel{\gamma}{\ra}& U\end{matrix}\]
commutes. We recall $U = \pi_* \sExt(\Odel, \Odel)$ and use
again the adjunction 
\[ \Hom(-, U) \cong \Hom(\pi^*(-),
\sExt(\Odel, \Odel))\] to compute the two sides of this square. Write
$\scrE$ for $\sExt(\Odel, \Odel))$ as a notational convenience.

The adjoint to the composition around the top is \[ \pi^*T[-1] \otimes
\pi^*T[-1] \stackrel{\pi^*(\gamma \otimes \gamma)}{\ra}
\pi^*(\pi_*\scrE \otimes \pi_*\scrE) \stackrel{[,]}{\ra} \scrE.\] The
right-hand map here, the commutator, is given in terms of the algebra
structure on $U$ which actually comes from a similar adjunction, so it
can be factorised
\[ \pi^*(\pi_*\scrE \otimes \pi_*\scrE) = \pi^*\pi_*\scrE \otimes
    \pi^*\pi_*\scrE \stackrel{p\otimes p}\ra\scrE \otimes \scrE \ra
    \scrE\] where $p$ is the adjunction unit $\pi^*\pi_*\scrE \ra
    \scrE$.  The composite $\pi^*T[-1] \otimes \pi^*T[-1] \ra
    \scrE\otimes \scrE$ obtained is by definition $\alpha \otimes
    \alpha$ so the whole map can be thought of as the commutator
    $[\alpha, \alpha]$.

The lower side of the square is adjoint to 
\[ \pi^*T[-1] \otimes \pi^*T[-1] \stackrel{\pi^*(\alpha_T)}{\ra} 
 \pi^*T[-1] \stackrel{\alpha}{\ra} \scrE.\]
Now the equality of these two compositions is just the STU identity
for the universal Atiyah class, proved in the previous section.

\noindent{\em Step 5}. 
The object $U$ acts on objects as follows.  We define a morphism
\[ \Odel \otimes \pi^*U \ra \Odel \]
in $D(X \times X)$ by taking the composition
\[ \Odel \otimes \pi^*\pi_*\scrE \ra \Odel \otimes \scrE \ra
\Odel\] using the unit of the adjunction and the natural
multiplication action of $\scrE$ on $\Odel$. This morphism
induces a natural transformation $- \otimes U \ra -$ which makes the
algebra $U$ act on $D(X)$. 

This is compatible with the action of $L=T[-1]$ on objects. To see
this we need to show that the diagram
\[ \begin{matrix} \Odel \otimes \pi^*T[-1]
  &\stackrel{\alpha}{\ra} & \Odel\\
\id \otimes \pi^*\gamma \downarrow&& \downarrow\\
\Odel \otimes U  &\ra & \Odel\end{matrix}\]
commutes. But equivalently we can transfer the $\Odel$s to the
other side and look at
\[ \begin{matrix} \pi^*T[-1]
  &\stackrel{\alpha}{\ra} & \scrE\\ \pi^*\gamma \downarrow&&
\downarrow\\ \pi^*\pi_*\scrE &\ra & \scrE\end{matrix}\] whose
commutativity is in fact the definition of $\gamma$.

\subsection{The PBW isomorphism}

We can finally construct the PBW isomorphism. Start with the canonical
map $\gamma: L \ra U$ coming from the Atiyah class. By tensoring it up
in the normal way it extends to an algebra homomorphism from the
tensor algebra $T(L)$ to $U$, and by composing with the
(non-algebra-morphism) symmetrisation map $S \ra T(L)$ we get our PBW
map.

To prove that this is an isomorphism of objects of $D(X)$, we can work
locally: such objects are just complexes of sheaves, and a map is an
isomorphism in $D(X)$ if it induces an isomorphism of cohomology {\em
sheaves}. Isomorphisms of sheaves can of course be checked locally in
an affine patch of $X$. Let $i: Y \hra X$ be an affine chart: from the
above remarks about locality of the Atiyah class (or by an abstract
functorial diagram-chase), we see that restricting $\PBW_X: S_X \ra
U_X $ gives the corresponding morphism $\PBW_Y: S_Y \ra U_Y$. So it is
enough to show that the PBW morphism is an isomorphism when $X$ is
affine.

To do this it helps to transfer from the category of coherent sheaves
on $X$ to the equivalent category of (left) $A$-modules, where
$A=\Gamma(\OO_X)$. Of course $\OO_X$ becomes the left regular module
$A$, $\Omega_X$ becomes the module of K\"ahler differentials
$\Omega^1_{A}$ and $T_X$ becomes the module of derivations
$\Der(A,A)$. The object $S$ is therefore represented by the exterior
algebra $\bigwedge_A \Der(A,A)$. (All tensor products in this section
are over $\C$ unless otherwise noted.)

Extending this dictionary, sheaves on $X \times X$ become
$A-A$-bimodules, that is $A^e$-modules, where $A^e$ is the enveloping
algebra $A \otimes A^{\op}$. In particular we have that $\OO_{X \times
X}$ corresponds to $A \otimes A$ (the free $A-A$-bimodule of rank $1$)
whereas $\Odel$ corresponds to $A$ as a bimodule. The cotangent
sheaf $\Omega_{X \times X}$ corresponds to the bimodule
\[ \Omega^1_{A^e} \cong \Omega^1_{A}\otimes A  \oplus A \otimes
\Omega^1_{A},\] where the two right-hand terms are of course the
pullbacks $\pi_1^*\Omega_X$ and $\pi_2^*\Omega_X$. Taking pushforward
$(\pi_1)_*$ simply corresponds to forgetting the right $A$-module action
of a bimodule, making it just a left $A$-module.

Computing the object $\sExt(\Odel, \Odel)$ is equivalent to computing
the bimodule $\sExt(A, A)$, where the $\sExt$ here is the derived
functor of internal hom in the category of $A-A$-bimodules'' (since
$A$ is commutative, the set of bimodule homomorphisms is itself a
bimodule). To compute it we can use a a resolution of the first factor
by free bimodules, such as the Hochschild (bar) complex:
\[ B(A) = \qquad \ra A^{\otimes n} \ra A^{\otimes n-1} \ra \cdots \ra
A\otimes A.\] Here the $A^{\otimes n}$ term is taken to be in degree
$2-n$, so we can in fact write \[B^{-n}(A) = A \otimes A^{\otimes n}
\otimes A \qquad n\geq 0.\] (Again, the tensor products are over $\C$
and the action of $A^e$ is on the outer factors). The differentials
are the usual Hochschild differentials:
\[ \partial (a_0 \otimes (a_1 \otimes \cdots \otimes a_n) \otimes a_{n+1}) =
\sum_{i=0}^n (-1)^ia_0 \otimes \cdots \otimes a_i.a_{i+1} \otimes \cdots
 \otimes a_{n+1}.\] The multiplication map $B^0(A) = A \otimes A \ra A$
 gives the resolving quasi-isomorphism $B(A) \ra A$.

The object $\sExt(A, A)$ is therefore represented by the complex of
bimodules $\Hom(B(A), A)$, but in order to see the algebra structure
most naturally we should resolve the second factor too, taking the
quasi-isomorphic complex $\Hom(B(A), B(A))$. We are really interested
in the object $U = \pi_*\sExt(A, A)$, which corresponds to the complex
of {\em left} $A$-modules $\Hom(B(A), B(A))$ (we just forget the right
module structure).

Now we can calculate explicitly the canonical map $\gamma: T[-1] \ra
\pi_*\sExt(A, A)$, which is the adjoint of the universal Atiyah class
map $\alpha$.  For this, we use Markarian's lemma 1 (from section 5). 

Since the Hochschild complex corresponds to a complex of trivial
sheaves on $X \times X$ (remember that $X$ is still assumed affine),
each of them has a global flat connection, namely the trivial
connection. $\nabla = d$. In the world of $A-A$-bimodules these
connections are given by the maps
\[ \nabla: B^{-n}(A) \ra B^{-n}(A) \otimes_{A^e} \Omega^1_{A^e} \]
which, using the decomposition of the module $\Omega^1_{A^e}$, becomes 
\[ \nabla: A \otimes A^{\otimes n} \otimes A \ra (\Omega^1_A
\otimes A^{\otimes n} \otimes A) \oplus (A \otimes A^{\otimes n} \otimes
\Omega^1_A)\] and, viewing the LHS as a free
$A^e$-module with basis $A^{\otimes n}$, is given explicitly by
\[ \nabla(a_0 \otimes (a_1 \otimes \cdots \otimes a_n) \otimes
a_{n+1}) = da_0 \otimes (a_1 \otimes \cdots \otimes a_n) \otimes a_{n+1}
+ a_0 \otimes (a_1 \otimes \cdots \otimes a_n) \otimes da_{n+1}.\]

Now the Atiyah class of $\Odel$ is represented by $\partial
\nabla - \nabla \partial$, where $\partial$ represents the Hochschild
differential, and we compute explicitly the difference
\[ (\partial \nabla - \nabla \partial) (a_0 \otimes (a_1 \otimes \cdots
  \otimes a_n) \otimes a_{n+1}).\] It's easy to see that the terms
coming from the Hochschild differential when $1 \leq i \leq n-1$
cancel out, leaving only the outer ($i=0,n$) terms, and we get the
answer
\[ a_0. da_1 \otimes a_2 \otimes \cdots \otimes a_n \otimes
a_{n+1} + (-1)^n a_0 \otimes a_1 \otimes \cdots \otimes da_n. a_{n+1}.\]

Since the actual universal Atiyah class $\alpha$ is the part involving
$\pi_1^*T$ we get a representation of \[ \alpha \in \Hom_{A^e} (B(A),
B(A) \otimes_{A^e} (\Omega^1_A \otimes A[1]))\] given by
\[ \alpha(a_0 \otimes (a_1 \otimes \cdots \otimes a_n) \otimes
a_{n+1}) = a_0. da_1 \otimes (a_2 \otimes \cdots \otimes a_n) \otimes
a_{n+1}.\]

Applying the adjunction we see that the map $\gamma: T \ra U[1]$ is
represented by the left module morphism
\[ \Der(A,A) \ra \Hom^1(B(A), B(A))\]
given by
\[ \gamma(\xi)(a_0 \otimes (a_1 \otimes \cdots \otimes a_n) \otimes
a_{n+1}) = a_0. \xi (a_1) \otimes (a_2 \otimes \cdots \otimes a_n) \otimes
a_{n+1}.\]

It follows that $\gamma^{\otimes n}: T^{\otimes n} \ra \sExt(\Odel,
\Odel)[n]$ is represented by a morphism
\[ \bigotimes_A^n\Der(A,A) \ra \Hom^n(B(A), B(A))\]
such that  
\[ \gamma^{\otimes n} (\xi_1 \otimes \cdots \otimes \xi_n) (a_0 \otimes (a_1 \otimes \cdots  
\otimes a_n) \otimes a_{n+1}) = a_0. \xi_1(a_1). \xi_2(a_2). \ldots
\xi_n(a_n) \otimes a_{n+1}.\] If we compose with the resolving
quasi-isomorphism $\Hom(B(A), B(A)) \cong \Hom(B(A), A)$ given by
composition on $B^0(A)$ with the multiplication map, we see that
$\gamma^{\otimes n}(\xi_1 \otimes \cdots \otimes \xi_n)$ lies in
$\Hom(B^{-n}(A), A)$ and is given by
\[  a_0 \otimes (a_1 \otimes \cdots  
\otimes a_n) \otimes a_{n+1} \mapsto a_0. \xi_1(a_1). \xi_2(a_2). \ldots
\xi_n(a_n). a_{n+1}\]

Finally we symmetrise over the $X\xi_i$ to obtain the map representing
the degree $n$ part of $\PBW$, 
\[ \textstyle\bigwedge^nT[-n] \ra \pi_* \sExt(\Odel, \Odel).\] But 
$\Hom_{A^e}(B(A), A)$ is the $n$th Hochschild cochain group, and the
symmetrised map we obtain is just the standard
Hochschild-Kostant-Rosenberg map which defines an isomorphism on cohomology
\[ \HKR: \textstyle\bigwedge^n_A \Der(A,A)) \ra HH^n(A,A).\]
This ends the proof that the PBW map is an isomorphism.

\subsection{The universal property of $U$}

To complete the proof of step 3, we need the following theorem of
Hinich and Vaintrob \cite{HV}.

\begin{thm}
Let $\scrC$ be a linear tensor category admitting infinite direct sums and
symmetrisers. Let $L$ be a Lie algebra in $\scrC$. Then there exists
a universal enveloping algebra $L \ra U(L)$ of $L$ in
$\scrC$. Furthermore there is a PBW isomorphism $S \cong U(L)$.
\end{thm}

Their proof of this theorem is essentially to start from the symmetric
algebra $S(L)$, which exists given the conditions on $\scrC$, and
then to redefine its product using a kind of universal algebraic
construction (and the language of operads). See also Deligne and
Morgan \cite{DMorg}.

Assuming the properties already proved in steps 1,2,4 and this
theorem, we can now complete the proof of the universal property of
our object $U$. By the Hinich-Vaintrob theorem, we know that in $D(X)$
a universal enveloping algebra $U(L)$ (with the desired universal
property) does exist. (We don't need to worry about infinite direct
sums: the symmetric algebra in our case is a finite sum.) All we need
to do is prove that our object $U$ is isomorphic, as an algebra, to
the Hinich-Vaintrob object $U(L)$. This is done by exploiting the
universal property of $U(L)$ as follows.

Since $U$ is an associative algebra and $L \ra U$ is a Lie algebra
homomorphism, this map extends to an algebra homomorphism $U(L) \ra
U$. If we can prove that this is an isomorphism of objects in $D(X)$
then we are done.

We have also the natural algebra homomorphism $T(L) \ra U(L)$ obtained
by extending $L \ra U(L)$ to a map of algebras, and the symmetrisation
morphism (viewing the symmetric algebra as a subspace of the tensor
algebra) $S(L) \ra T(L)$. Consider the composition of these with our
map $U(L) \ra U$:
\[ S(L) \ra T(L) \ra U(L) \ra U.\]
The composite of the first two maps is the universal PBW isomorphism
$S \cong U(L)$ constructed by Hinich and Vaintrob. On the other hand,
the composite of the latter two morphisms is the natural map $T(L) \ra
U$ extending $L \ra U$ to a map of algebras, and thus the whole
composition is by definition {\em our} PBW isomorphism $S \ra
U$. Therefore the final map $U(L) \ra U$ is also an isomorphism.

\subsection{Invariant parts}
The PBW isomorphism between $S$ and $U$ restricts to their invariant
parts. In standard Lie theory the invariant part of a module can be
thought of as $V^\g\cong \Hom_\g(\C,V)$, and this gives the right way
to generalise the notion to the categorical setting: the invariant
part of an object $A \in D(X)$ is $\Hom_{D(X)}(\OO_X, A)$, which is a
cohomology space.

In our context we see that 
\[ \Hom(\OO_X, S) =
H^\gradestar(\textstyle\bigwedge^\gradestar T)\] is the cohomology of
polyvector fields on $X$, called $\HT^*(X)$ by Kontsevich. The degree
$k$ piece is
\[ HT^k(X)=\bigoplus_{i+j=k}H^i(\textstyle\bigwedge^j T).\]
It is also worth identifying the invariant part of the symmetric
algebra of the {\em dual} $\Omega[1]$ of the Lie algebra $T[-1]$,
which is the usual Dolbeault cohomology of $X$ but with a grading
shift: in this context its natural part of degree $k$ is
\[H^0(\bigoplus_j
\textstyle\bigwedge^jT^\dualstar [k+j]) =
\bigoplus_{i-j=k}H^i(\textstyle\bigwedge^jT^\dualstar ).\] There is an
obvious ``cap product'' action of this cohomology ring on
$\HT^\gradestar(X)$.

On the other hand, the invariant part of $U$ is
\[ \HH^\gradestar(X)= \Ext_{X \times X}^\gradestar
(\Odel, \Odel),\] the $\Ext$-algebra of the structure sheaf of the
diagonal in $X \times X$, with the Yoneda product as algebra
structure. This should be thought of as the Hochschild cohomology of
the manifold $X$: the usual definition of the Hochschild cohomology of
an algebra $A$ is as $\Ext^\gradestar _{A\otimes A^{\text{op}}}(A,A)$
--- that is, the $\Ext$-algebra in the category of $A-A$-bimodules,
and the above definition of $\HH^*(X)$ is clearly the sheaf-theoretic
analogue.

The PBW isomorphism between $S$ and $U$ induces an isomorphism $\HKR:
\HT^*(X) \cong \HH^*(X)$. This version of the
Hochschild-Kostant-Rosenberg theorem is originally due to Gerstenhaber
and Schack \cite{GS}. Kontsevich showed how to alter it into an {\em
algebra} isomorphism, and we discuss this in section 8. 

\begin{remk} Using Hinich and Vaintrob's results \cite{HV}, the Hochschild
cohomology $\HH^\gradestar(X)$ can also be described ``externally'' as
the quotient of $\bigoplus_{i+j=n}H^i(X,T^{\otimes j})$ by relations
saying that the action of the Atiyah class equals the commutator
(i.e.\ the relations for a universal enveloping algebra).
\end{remk}

\subsection{Alternative approach to $U$}

There is a slightly different way to define $U$ using the
Grothendieck-Verdier functor. In some ways this is more natural, but
it looks even more abstract.

Recall that $\Delta\colon X\to X\times X$ denotes the diagonal
embedding and $\pi\colon X\times X \to X$ denotes projection onto the
first factor. Then we can set $U=\Delta^!\Odel$. (Recall from
Section~\ref{Section:SheafOperations} that $\Delta^!$ is the right
adjoint of $\Delta_*$.) Then we can see that there are isomorphisms of
functors:
\begin{align*}
  \Delta_*\cong \Odel\otimes \pi^*(-);\qquad
  \Delta^!\cong \pi_*\sExt(\Odel,-).
\end{align*}
This follows from the projection formula and the fact that $\pi\circ
\Delta=\id$. We have
  $$\Odel\otimes \pi^*(-)
      \cong\Delta_*(\OO_X\otimes\Delta^*\pi^*(-))
      \cong \Delta_*(-).$$
So we can write $\Delta_*$ as the composition
$(\Odel\otimes-)\circ \pi^*$.  By the uniqueness of adjoints
we can write the right adjoint $\Delta^!$ as the composite of the
right adjoints of the components, viz:
$$\Delta^!\cong\pi_*\circ(\sExt(\Odel,-))=\pi_*\sExt(\Odel,-).$$
Note that the adjunctions, such as $\Delta_*\Delta^!\to \id$,
translate into the composition of adjunctions, such as
$\Odel\otimes
\pi^*\pi_*\sExt(\Odel,-)\to\Odel\otimes
\sExt(\Odel,-)\to\id$.

In this approach, to give the action of $U$ via a natural
transformation $-\otimes U\to \id$, it suffices to give a map
$\Delta_*U\to \Odel$.  As $U=\Delta^!\Odel$ we take the adjunction
$\eta\colon\Delta_*\Delta^!\Odel\to \Odel$. If $U$ is thought of as
$\pi_*\sExt(\Odel,\Odel)$ then the map is the composition
$\Odel\otimes \pi^*\pi_*\sExt(\Odel,-) \to\Odel\otimes
\sExt(\Odel,-)\to\id$ of the two basic adjunctions.

Similarly, the canonical Lie algebra homomorphism $T[-1]\to U$ can be
defined as the righht adjoint of the universal Atiyah class morphism
$\alpha\colon \Delta_*T[-1]\to \Odel$. 

\begin{remk} The associative algebra object $U$ which acts
on the objects of $D(X)$ has been constructed from functors on derived
categories induced by the diagonal map $\Delta\colon X\to X\times X$
and the projection map $\pi\colon X\times X \to X$. Starting with a
finite group $G$ an analogous construction can be performed using
functors on representation categories induced by the diagonal map
$\Delta\colon G\to G\times G$ and the projection map $\pi\colon
G\times G \to G$, in this case the resulting algebra object in the
representation category of $G$ which acts on everything in the
category is nothing other than the group algebra of $G$, equipped with
the adjoint action.  Details of this will appear elsewhere.
\end{remk}




\pg

\section{Further weight systems}

We now look at the roles played in Vassiliev theory by the symmetric
and universal enveloping algebras of a Lie algebra, and construct
weight systems from complex symplectic manifolds in this context.
We begin by introducing a new space of diagrams resembling $\A$. 

Define $\B$ to be the vector space spanned by not-necessarily
connected unitrivalent graphs with the same vertex-orientation
convention at their trivalent vertices, and subject to the
antisymmetry and IHX relations as before.
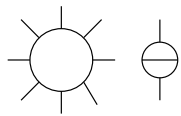
\[\vpic{Belement}\]
Again we use the total number of trivalent and univalent
vertices as grading, though we can also {\em bigrade} the algebra
and write $\B^{v,l}$ for the part with $v$ internal trivalent
vertices, and $l$ legs.  The vector space $\B$ is naturally a
commutative algebra via $\amalg$, the disjoint union of diagrams.

There is an isomorphism of graded, complex vector spaces
\[\chi\colon \B \ra \A\]
given by taking an $l$-legged diagram in
$\B$ to the average of the $l!$ diagrams obtained by attaching its
legs in all possible orders to an oriented circle (see \cite{BN}).
The isomorphism $\chi$ is not an algebra isomorphism, so it is
sometimes convenient to regard $\B$ and $\A$ as one space, using
$\chi$, which has two competing products. However there {\em is}
an interesting algebra isomorphism between $\A$ and $\B$ which is
described below.


\subsection{Further weight systems from Lie algebras}
\label{Subsection:FurtherLAweightsystems} We can now construct
further weight systems, and will encapsulate them all in the
following theorem.  We give a proof in this familiar context, as
this proof will go over pretty much exactly to the complex
symplectic context in
Section~\ref{Subsection:FurtherHSMweightsystems}.
\begin{thm}[See \cite{BN}]
\label{Thm:FurtherLAWeightSystems}
Suppose that $\g$ is a finite-dimensional metric Lie algebra. Let
$S(\g)$ be its symmetric algebra, and $U(\g)$ be its universal
enveloping algebra.
\begin{enumerate}
\item There is a graded, multiplicative weight system
$w_\g\colon\B\to S(\g)^\g$.
\item There is a graded, multiplicative weight system
$w^\g\colon\B\to S(\g^\dualstar)^\g$.
  \item Given $V$ a finite-dimensional representation of $\g$, there
  is a multiplicative weight system $w_V\colon \A\to \End(V)$;
  composing with the trace we get a weight system $w_V\colon \A\to
  \C$.  %
\item There is a multiplicative weight system $w_\g\colon\A\to
U(\g)^\g$.  If $V$ a finite-dimensional representation of $\g$
then composing with the natural map $U(\g)^\g\to \End(V)$ gives
the weight system in {\rm (3)} above.
\item The maps $\chi$ and $\PBW$ correspond in the sense that the
 following diagram commutes:
\[
 \begin{array}{ccc}
  \B  & \stackrel{\chi}{\longrightarrow} &\A\\
  \llap{$w_\g$}\downarrow  && \downarrow \rlap{$w_\g$}\\
  S(\g)^\g &  \stackrel{\PBW}{\longrightarrow}&U(\g)^\g
 \end{array}
\]
\end{enumerate}
\end{thm}
\begin{proof}
These weight systems are all defined essentially by taking any
Morse, planar projection of a representing graph and viewing it
as a morphism in the category of $\g$--modules.  That it will be
independent of the choice of projection and morsification
 is due precisely to the axioms of a Lie algebra object in a category and
of modules over it.  This is the work of Vogel and Vaintrob.  The work here is
 really in properly identifying the target.

Parts (1) and (2) are straightforward. Given a diagram in $\B$,
represent it in the plane with all legs pointing upwards (in case (2),
point them downwards and make obvious alterations). The legs will have
to be ordered arbitrarily from left to right to do this. The picture
defines an element of $\Mor_\g(\C,\TT(\g))$; composing with the
canonical map $\TT(\g) \ra S(\g)$ gives a result independent of leg
ordering. Multiplicativity in the first case follows by placing
diagrams side-by-side. In the second case, the target space can be
thought of as the algebra $\Mor_\g(\C, S(\g^\dualstar ))$ and the map
is mutliplicative, but this is less important.

For part (3), cut the diagram in $\A$ at some point of its
oriented circle, and open it out to an upward-oriented interval,
with attached graph drawn to the right. This picture defines an
element of $\Mor_\g(V,V)$. The result is independent of the location
of the cut by the standard argument from Bar-Natan, which we draw
here.
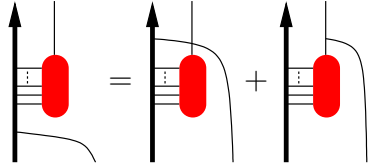
\[\vpic{cyc1} = \vpic{cyc2} + \vpic{cyc3}
\]
This expresses the fact that the ``oval with
legs'' (representing any graph with legs)  is an
$L$-module map.
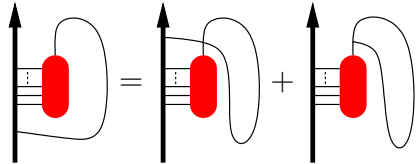
\[\vpic{cyc4}=\vpic{cyc5}+\vpic{cyc6}
\]
This follows by applying the Casimir and
metric, and untangling the pictures suitably. Note that although
$\Mor_\g(V,V)$ need not be a commutative algebra, the image of the
weight system is a commutative subalgebra.

Part (4) is only a little different. This time we cut the circle
and draw the remaining interval {\em horizontally}, pointing to
the right, with the rest of the graph below it. Removing the
oriented interval gives a graph with legs ordered from left to
right. This defines an element of $\TT(\g)^\g$ which projects to
something in $U(\g)^\g$. The IHX relations in $\A$ are clearly
respected, and the STU relations also because of the universal
property of the canonical morphism $\g \ra U(\g)$. Independence of
the point of cutting follows from the same pictorial argument as
above.

The comparison with part (3) arises as follows.
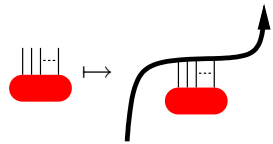
\[\vpic{uea1}\mapsto\vpic{uea2}\]
Part (5) is now a straightforward check.
\end{proof}

The above construction works for any metric Lie algebra object in a
category, so the case of complex symplectic manifolds will follow
naturally. One point worth making is that it is clear from this
construction that the notion of `invariant part' of a module $M$
should be the hom-set $\Hom_L(\one, M)$.


\subsection{The Duflo isomorphism and wheeling}

We described earlier the PBW isomorphism between spaces of invariants
$S(\g)^\g$ and $U(\g)^\g$. Each of these spaces is a commutative
algebra, but the $\PBW$ map is not generally an algebra
isomorphism. There {\em is} however an algebra isomorphism, the {\em
Duflo isomorphism}, between $S(\g)^\g$ and $U(\g)^\g$ --- for
semisimple Lie algebras it is equivalent to Harish-Chandra's
isomorphism, but in the more general form it is due to Duflo
\cite{Duflo}.
 
To define it, consider the invariant polynomial function
$s_i(x)=\tr(\ad(x)^i)$ on $\g$ as an element of the dual symmetric
algebra $S^i(\g^\dualstar)^\g$. This space acts on $S(\g)^\g$ by the
symmetrised contraction map.  We will think of it as a kind of cap
product, and will write $f\cap -: S^\gradestar (\g)^\g \ra
S^{\gradestar -i}(\g)^\g$ for the action of $f\in S^i(\g^\dualstar
)^\g$.

Define the {\em modified Bernoulli numbers\/}
$\{b_{2i}\}_{i=1}^\infty$ by the power series
\[   \sum_{i=1}^\infty b_{2i}x^{2i}=
\frac12\log\frac{\sinh(x/2)}{x/2}, \] and define the Duflo power
series
\[ j^{\frac12}=\exp\sum b_{2i}s_{2i}\]
in the completion of $S(\g^\dualstar)^\g$.

\begin{remk}
This function plays a very important role in the Weyl character
formula, amongst other things. For a semisimple Lie algebra we could
identify the invariant polynomials $S^i(\g^\dualstar )^\g$ with
$H^{2i}(BG)$, so that $s_i$ would correspond to $i!$ times the $i$th
term of the Chern character of the vector bundle on $BG$ corresponding
to the adjoint representation. In $H^\gradestar (BG)$, it corresponds
to the equivariant $\hat A$-genus of the complex adjoint
representation \cite{BGV}.
\end{remk}

The Wheeling Theorem of Bar-Natan, Le and Thurston \cite{BLT} is
a strange and deep property of the algebras $\A$ and $\B$ which
corresponds to the Duflo isomorphism.

The algebra $\B$ acts on itself by a leg-gluing operation which we
will here denote by a ``cap'' notation.  Thurston uses a ``hat'' or
``differential operator'' notation.  This operation is defined on
diagrams $C$ and $D$ by
\[  C \cap D =\sum \text{all ways of joining all of the legs of $C$ to some of the legs of $D$.}\]
If $C$ has more legs than $D$ then, $C \cap D$ is zero. The
capping operation $ - \cap - : \B\otimes\B\to\B$ is not a graded
map, but if $\B$ is given an alternative ``Euler characteristic''
grading, namely $\B_n=\bigoplus_{l-v=n}\B^{v,l}$, then $ - \cap -:
\B_\gradestar \otimes\B^\gradestar \to\B^\gradestar $ is graded.

Let $w_l$ denote the wheel with $l$ legs and let $\Omega\in \B$ be the
{\em wheeling element} given by the following
formula.
\[  \Omega :=\exp_{\amalg} \sum_{i=1}^\infty b_{2i}w_{2i}\in \B. \]
It is in the subspace $\B_0$, so that the {\em wheeling map}
$\Omega\cap -: \B^\gradestar \ra \B^\gradestar $ is a graded map.
Note that although $\Omega$ really lives in the completion of
$\B$, there is no ``convergence problem'' when we define the
wheeling map.

\begin{thm}[Wheeling Theorem \cite{BLT}]
The composition of the wheeling map   $\Omega\cap {-}$  with the
symmetrisation map is an algebra isomorphism $\B\to\A$.
\end{thm}

If $\g$ is a finite-dimensional metric Lie algebra then we can combine
the weight systems and the above isomorphisms into the following
commutative diagram, whose top and bottom rows are both algebra
isomorphisms. (Note in particular that $\Omega$ maps to $j^{1/2}$, a
fact originally pointed out in \cite{BGRT}.)

\[\begin{array}{ccccc}
\B
  & \stackrel{\Omega\cap-}{\longrightarrow}
  & \B
  &   \stackrel{\chi}{\longrightarrow}
  & \A\\
\downarrow & & \downarrow & & \downarrow\\
S(\g)^\g
  &  \stackrel{{j^{\frac12}}\cap-}{\longrightarrow}
  &  S(\g)^\g
  & \stackrel{\PBW}{\longrightarrow}
  &  U(\g)^\g
\end{array}\]


Note that this whole diagrammatic method could be used profitably for
handling higher graph cohomology and $L_\infty$-algebras
--- Kapranov uses the language of operads. However, for our present
purposes, this extra structure is not important.

\subsection{Further weight systems from complex symplectic manifolds}
\label{Subsection:FurtherHSMweightsystems}

We can now translate Theorem~\ref{Thm:FurtherLAWeightSystems} directly
into the context of complex symplectic manifolds.  Thus the
following theorem has a parallel proof to that of
Theorem~\ref{Thm:FurtherLAWeightSystems} using the structures
discussed in the above subsection.

\begin{thm}
Suppose $(X,\omega)$ is a complex symplectic manifold.
\begin{enumerate}
\item There is a bigraded, multiplicative weight system
$\RW_X\colon  \B^{\gradestar,\gradestar}
\ra H^\gradestar (\bigwedge^\gradestar T)$.

\item There is a bigraded, multiplicative weight system
$\RW^X\colon  \B^{\gradestar,\gradestar}
 \ra
H^\gradestar (\bigwedge^\gradestar \Omega )$.

\item Given an object $A \in \tilde D(X)$, there is a graded
multiplicative weight system $\A \ra \Ext^\gradestar (A,A)$, and
composing with the trace we get a weight system $\A\to H^\gradestar
(\OO_X)$.

\item There is a graded multiplicative weight system $\A \ra
\HH^\gradestar(X) = \Ext^\gradestar _{D(X \times X)}(\Odel,
\Odel)$. If $A \in D(X)$ and we compose with
the natural map $\HH^\gradestar(X)  \ra \Hom_{D(X)}(A, A)$, we recapture
the weight system from {\rm (}3{\rm )}.

\item The HKR map $\HT^\gradestar(X)  \ra \HH^\gradestar(X) $ induces the
following commutative diagram of vector spaces:
\[\begin{array}{ccc}
\B&\stackrel\chi\longrightarrow&\A\\
\downarrow&&\downarrow\\
\HT^\gradestar (X)&\stackrel{HKR}\longrightarrow&\HH^\gradestar(X)
\end{array}.\]
\end{enumerate}
\end{thm}


\subsection{Wheels and wheeling for complex symplectic manifolds}

The wheeling theorem for complex symplectic manifolds takes the
following form:

\begin{thm}
Let $X$ be a complex symplectic manifold. Then there is a
commutative diagram
\[\begin{array}{ccccc}
\B
  & \stackrel{\Omega\cap-}{\longrightarrow}
  & \B
  &   \stackrel{\chi}{\longrightarrow}
  & \A\\
\downarrow & & \downarrow & & \downarrow\\
\HT^\gradestar(X)
  &  \stackrel{{\hat A^{\frac12}}\cap-}{\longrightarrow}
  &  \HT^\gradestar(X)
  & \stackrel{\text{HKR}}{\longrightarrow}
  &  \HH^\gradestar(X)
\end{array}\]
in which the two rows are algebra isomorphisms.
\end{thm}
\begin{proof}
The fact that the bottom line is an algebra isomorphism is due to
Kontsevich \cite{KontsDefQuant1}.  Note that it holds for any complex
manifold.

The intertwining weight system maps are only defined when the manifold
is complex symplectic. That the square on the left commutes
follows from the following lemma, independently computed by Hitchin
and Sawon \cite{Hitchin-Sawon}.
\end{proof}

\begin{lemma}
$RW^X(\Omega) = \hat A^{1/2}(TX) \in \bigoplus H^{2k}(\bigwedge^{2k}T^\dualstar )$.
\end{lemma}
\begin{proof}
All we need is that an $l$-leg wheel $w_l$, with its hub oriented and
labelled by a locally-free sheaf $E$ and legs pointed downwards maps
under the weight system to
\[ \RW^{X,E}(w_l)= \tr(\tilde F_E^{l}) \in H^l(\textstyle\bigwedge^lT^\dualstar ).\]
This is a restated lemma of Atiyah \cite{A}: we are using the
Dolbeault point of view, in which $\tilde F_E \in
\Omega^{1,1}(\End(E))$ is the renormalised curvature form $\tilde F_E
= -1/2\pi i F_E$ of a smooth hermitian connection on $E$.

In particular, wheels in the honest algebra $\B$ correspond to the
above case when $E=T[-1]$. The only effect of the degree shift here is
to make the trace negative (it is really the supertrace of an odd
object) and thus the $l$-wheel gives $-l!  \ch(TX)$. (We did not
specify an orientation on the hub of the wheel, but we can introduce
one arbitrarily when comparing the definitions: odd wheels in $\B$ are
zero, corresponding to the vanishing of the odd Chern classes of a
complex symplectic manifold.)

Now recall that $\Omega \in \B$ is defined as
\[ \Omega = \exp \sum b_{2n}w_{2n} \]
with
\[ \sum b_{2n}x^{2n}= \frac12\log \frac{\sinh(x/2)}{x/2}.\]
Under the weight system, each $2k$-wheel goes to a term $-\tr
(\tilde F^{2k})$, where we take trace in the fundamental
representation. Because the weight system is multiplicative,
disjoint union becomes cup product: we get $\exp \tr -\frac12 \log
\frac{\sinh(\tilde F/2)}{\tilde F/2}$, which is just the
Chern-Weil definition of the $\hat A^{1/2}$.
\end{proof}

It is worth pointing out that in general, a wheel whose hub is
labelled by an object $A \in D(X)$ maps to the characteristic class
$l! \ch_l(A)$, by which we mean the alternating sum of the terms
$\ch_l(\scrE)$ for the sheaves in the complex $A$, multiplied by
$l!$. Evaluations of this nature will feature in our future work but
are not needed for now. 

Note that for a complex symplectic manifold, the $\hat A$ and Todd
classes are equal, because $T \cong T^\dualstar $ means that the line
bundles appearing in the splitting principle occur in conjugate pairs,
and the first Chern class is therefore zero. But the $\hat A$ class is
the correct one to use, as it appears in Kontsevich's theorem, which
holds for {\em any} complex manifold (even if $c_1$ is not zero).

The appearance of the class $\hat A^{\frac12}$ requires further
investigation. Does it have a meaning in index theory: for
example, is there an interesting class of manifolds whose $\hat
A^{\frac12}$-genus is integral?   According to Sawon
\cite{Sawonthesis}, it is not integral for compact hyperk\"ahler
manifolds.

\pg

\section{Ribbon categories and link invariants}

In this section we combine the complex symplectic manifold weight
systems with the Kontsevich integral to obtain an interesting {
ribbon category}, from which link invariants may be obtained by
the standard methods of Turaev.

\subsection{Ribbon categories}

In Section~\ref{Subsection:SymmetricTensorCats} we introduced
symmetric tensor categories; braided tensor categories are to the
braid groups as symmetric tensor categories are to the symmetric
groups.  A tensor category is\/ {\em braided} if there is a
natural isomorphism $\tau$, the {braiding}, between $\otimes$ and
$\otimes \circ \sigma$ (where $\sigma$ is the obvious flip functor
$\scrC \times \scrC \ra \scrC \times \scrC$), satisfying the
hexagon relation
\[
  \tau_{A, B \otimes C} =(\id_B \otimes \tau_{A,C}) \circ (\tau_{A,B}
  \otimes \id_C),
\]
where the component $A\otimes B\to B\otimes A$ of the natural
isomorphism is written $\tau_{A,B}$. The hexagon would be more
visible if we hadn't dropped the associators from the notation.
The braiding depicted as follows:
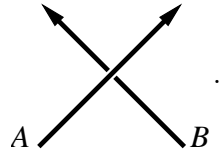
\[\vpic{braid}.\]
 Combining the hexagon
condition with naturality of $\tau$ yields the Yang-Baxter (or
braid, or Reidemeister III) equation
\[
  (\tau_{B,C}\otimes \id)\circ(\id\otimes\tau_{A,C})
  \circ(\tau_{A,B}\otimes\id) =(\id\otimes\tau_{A,B})
  \circ(\tau_{A,C}\otimes\id)\circ(\id\otimes\tau_{B,C}).
\]
The reader is strongly encouraged to draw the picture.  A
symmetric tensor category is a braided tensor category in which
the square of $\tau$ is the identity.

In a braided tensor category, there is an action of the $n$ string
braid group on the $n$th tensor power of any object. In a
symmetric tensor category, this factors through an action of the
symmetric group.

A {\em ribbon category} (or {\em balanced rigid braided tensor
category}) is a braided tensor category with a\/ {\em twist\/}
$\theta$ which is a natural automorphism of the identity functor
that commutes with duality and interacts with tensor product
according to the formula
\[
  \theta_{A \otimes B} = \tau_{B,A}\tau_{A,B} (\theta_A \otimes
  \theta_B).
\]
Using this it is possible to make natural isomorphisms $A \cong
A^{\dualstar \dualstar }$ and $A^\dualstar  \cong
{{}^\dualstar\!\! A}$ in ways compatible with tensor product, and
which can be neglected notationally.

The idea behind ribbon categories is that morphisms are thought of
as being two-sided ribbons, rather than strings.  The twist
$\theta$ represents a full-twist of the ribbon and is illustrated
diagrammatically as below left or sometimes as on the right as it
is easier to draw.
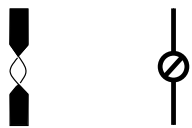
\[\vpic{band}\qquad \qquad \vpic{lollipop}\]
The reader is invited to discover the topological identity lurking
in the tensor product interaction described above. Essentially by
definition, a ribbon category gives rise to invariants of framed
links in the following way.  If the components of a framed link
are coloured with objects in the ribbon category then any morse
diagram of the link can be interpreted as a morphism from the unit
object to itself.  This element of $\Mor(\one,\one)$ is an
invariant of the coloured, framed, oriented link.  Quantum
invariants can be obtained in the this fashion using
representation categories of quantum groups.

\subsection{The ribbon structure of $\tilde D(X)$}

The construction of an interesting ribbon structure on $\tilde D(X)$
is parallel to the construction of an interesting braided structure on
the category of modules over a finite-dimensional metric Lie algebra
using the Knizhnik-Zamolodchikov equation. One starts with the usual
(symmetric) category $\gmod$, tensors with $\C[[\hbar]]$, and uses the
KZ equation to introduce a new, interesting braiding structure. The
resulting category turns out to be equivalent to the category of
representations of the quantum group $U_\hbar \g$.  See Bakalov and
Kirillov \cite{BK}, for example.

In our case, we start with the derived category $D(X)$ of a
complex symplectic manifold. We know this is a symmetric tensor
category under derived tensor product and the usual flip map, with
the structure sheaf as identity, and is a ribbon category when we
bring in the derived duals of objects and the natural contraction
maps.  We ``tensor with $\C[[\hbar]]$'', replacing $D(X)$ by the
extended version $\tilde D(X)$, in which the shift $[2]$ plays the
role of $\hbar$. Then we use the Kontsevich integral, in the
tangle-functor version of Le and Murakami, to intoduce the ribbon
structures. The result is that the derived category of coherent
sheaves on an complex symplectic manifold can be quantized in
exactly the same way as the category of representations of a
metric Lie algebra. This is our final main theorem.

\begin{thm}
For a complex symplectic manifold $X$, the extended derived
category $\tilde D(X)$ has a natural non-symmetric ribbon tensor
category structure.
\end{thm}
\begin{proof}
All we have to do is define the various structure morphisms and check
the identities. Explicitly, we need to compute the associator
$\Phi_{A,B,C}$, the braiding $\tau_{A,B}$, the duality morphisms
$\epsilon_A, \iota_A, \epsilon_A', \iota_A'$ and the twist
coefficients $\theta_A$.

According to Le and Murakami, the Kontsevich integral defines a
representation of the category of {\em non-associative tangles} (also
known as {\em quasi-} or {\em $q$-tangles}). This category is
generated by morphisms corresponding to exactly the things we need
above, and their explicit Kontsevich integrals can be found in Le and
Murakami. Each is a formal power series of diagrams consisting of
chords based on an underlying collection of oriented intervals: three
for $\Phi$, two for $\tau$, and one for the other morphisms.  For
example, $\Phi$ has an expression as a power series in the two
diagrams shown below, composed vertically and sometimes thought of as
non-commuting indeterminates:
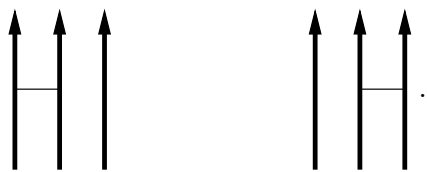
\[ \vpic{chordl} \qquad \qquad \qquad \vpic{chordr}.\]

Having obtained these power series, we label the vertical strings of
the diagram by objects of $\tilde D(X)$ and evaluate using the weight
systems. In particular, given two objects $A,B$ of $\tilde D(X)$, let
$H_{A,B} \in
\Mor_{\tilde D(X)} (A \otimes B, A \otimes B)$ and $C_A \in
\Mor_{\tilde D(X)}(A, A)$ be the morphisms
corresponding to the following graphs:
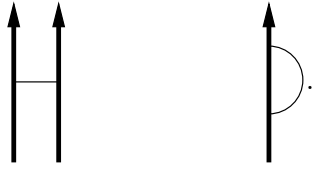
\[ \vpic{chordH} \qquad \qquad \qquad \vpic{chordcas}.\]

Each is really an element of $\Ext^2$, that is $H_{A,B} \in \Ext^2(A
\otimes B, A \otimes B)$ and $C_A \in \Ext^2(A,A)$. These ``chord''
and ``Casimir'' elements are really all we need to evaluate the
Kontsevich integrals. For example, the new braiding morphism
$\tau_{A,B}$ may be described as
\[ \tau_{A,B} = \tau_{\text{old}} \circ \exp(H_{A,B}/2) \in \Ext^\gradestar (A\otimes B, B \otimes A) =
\Mor_{\tilde D(X)}(A \otimes B, B \otimes A), \]
where $ \tau_{\text{old}}$ is the original symmetric braiding.  The
associator $\Phi_{A,B,C}$ is written as a polynomial in the
non-commuting variables $H_{A,B} \otimes \id_C$, $\id_A \otimes
H_{B,C}$. (Note that the power series become truncated because of the
boundedness of the $\Ext$-groups.) The other morphisms depend directly
on the Casimir. For example the framing twist is
\[ \theta_A = \exp(C_A/2) \in \Ext^\gradestar (A,A) = \Mor_{\tilde D(X)}(A,A).\]

The fact that all the relations of a ribbon category are satisfied is
then automatic from the topological invariance of the Kontsevich
integral of framed oriented tangles.
\end{proof}

A few remarks are now in order.  Firstly, one can multiply the
symplectic form by $\hbar$ and then ``take the limit $\hbar \ra
0$'' to recapture the original symmetric structure on $\tilde
D(X)$.

Secondly, we don't know whether it is possible to make a ``gauge
transformation'' (in the manner of Drinfeld) to a form which is
strictly associative but has a more complicated braiding (as in the
case of quantum groups). This theorem (discussed in \cite{BK}) does
not seem to have a purely geometrical formulation which can be carried
over into our context.

Thirdly, the above construction actually goes through for a
holomorphic Casimir manifold. In fact, only chord diagrams are
used in the construction, so a complex manifold $X$ with a
suitable ``$r$-matrix''  would be sufficient.

Finally, we observe that the braided structure is completely different
from the braid group actions on derived categories constructed by
Seidel and Thomas or  Rouquier.


\pg

\section{Conclusion}

In the table below we present a dictionary giving a translation
between the worlds of Chern-Simons theory (derived from usual Lie
algebras) and Rozansky-Witten theory (derived from complex symplectic
manifolds).

\begin{table}[h]
\renewcommand{\arraystretch}{1.5}
{\small
\begin{tabular}{|c||c|c|}
\hline
  & Chern-Simons
  & Rozansky-Witten
\\\hline\hline
\raggedright Category
  & Vector spaces
  & $D(X)$, derived category of $X$
\\\hline
\raggedright Lie algebra object
  & $\g$
  & $T[-1]$, shifted tangent bundle
\\\hline
\raggedright Modules
  &$ \rho\colon\g\to\operatorname{End}(V)$
  &objects $A$ of $D(X)$
\\\hline
\parbox[c][1.5\height]{9em}{Invariant part of\\ enveloping
algebra}
  & $Z(\g)\cong U(\g)^\g$
  & $\operatorname{Ext}^\gradestar _{X\times
          X}(\Odel, \Odel)$
\\\hline
\parbox[c][1.5\height]{9em}{Invariant part of\\ symmetric algebra}
  &  $S(\g)^\g$
  &   $H^\gradestar (X,\bigwedge^\gradestar  T) $
 \\\hline
\raggedright Wheeling theorem
  & \parbox{20ex}
{\center{$S(\g)^\g\stackrel\cong\longrightarrow
      Z(\g)$ \\ Duflo}}
  & \parbox[c][1.5\height]{38.5ex}{\center{$H^\gradestar (X,\bigwedge^\gradestar  T) \stackrel\cong\longrightarrow
    \operatorname{Ext}^\gradestar _{X\times X}(\Odel, \Odel)$ \\Kontsevich}}
\\\hline
\raggedright Invariant metric
  & $\g\otimes\g\to\mathbb C$
  & $\omega\in \Gamma(\bigwedge^2T^\dualstar)$
\\\hline
\raggedright Universal knot invariant
  & $\{\text{knots}\}\to Z(\g)[[h]]$
  & $\{\text{knots}\}\to\operatorname{Ext}^\gradestar _{X\times X}
                                         (\Odel,\Odel)$
\\\hline

\raggedright Knot invariant from module
  &$\{\text{knots}\}\to \mathbb C[[h]]$
  &$\{\text{knots}\}\to H^\gradestar (X,\mathcal O)$
\\\hline
\raggedright Ribbon category
  &$U_{\hbar}(\g)$-modules
  &$\tilde D(X)$
\\\hline

\end{tabular}}\bigskip

\end{table}

This table only goes as far as the knot invariants arising in each
theory, but we will extend this table into a correspondence between
the full TQFTs in a subsequent paper (a sketch of this appears in
\cite{R}).

There are many interesting questions arising from the existence of the
Rozansky-Witten invariants, their similarity with Chern-Simons
constructions, and their potential applications in knot theory. In
Roberts and Sawon \cite{RS} we mentioned many of these, so there is
little point repeating them here. 


\pg


\begin{thebibliography}{BGRT}

\bibitem[A]{A} M.~F.~Atiyah. {\em Complex analytic connections in
        fibre bundles}, Trans.\ Amer.\ Math.\ Soc.\ 85 (1957),
        181-207.

\bibitem[BK]{BK} B.~Bakalov, A.~Kirillov. {\em Lectures on tensor
        categories and modular functors} University Lecture Series 21,
        AMS (2001).

\bibitem[BN]{BN} D.~Bar-Natan. {\em On the
        Vassiliev knot invariants}, Topology 34 (1995), 423-472.

\bibitem[BGRT]{BGRT} D.~Bar-Natan, S.~Garoufalidis, L.~Rozansky,
        D.~P.~Thurston. {\em Wheels, wheeling, and the Kontsevich
        integral of the unknot}, Israel J.\ Math.\ 119 (2000),
        217-237.

\bibitem[BLT]{BLT} D. Bar-Natan, T.~Le, D.~Thurston, {\em Two
applications of elementary knot theory to Lie algebras and Vassiliev
invariants}. Geom. Topol. 7 (2003), 1-31.

\bibitem[BGV]{BGV} N.~Berline, E.~Getzler, M.~Vergne. {\em Heat
        kernels and Dirac operators}, Grundlehren der Mathematischen
        Wissenschaften 298, Springer (1992).

\bibitem[B]{Beauville} A.~Beauville. {\em Vari\'et\'es K\"ahleriennes
        dont la premi\`ere classe de Chern est nulle}, J.\
        Differential Geom.\ 18 (1983), 755-782.

\bibitem[Ca]{CaldararuII} A.~C\u ald\u araru. {\em The Mukai pairing,
         II: the Hochschild-Kostant-Rosenberg isomorphism}, Adv. Math.
         194 (2005), 34-66.

\bibitem[CP]{CP} V.~Chari, A.~Pressley. {\em A guide to quantum
        groups}, Cambridge University Press (1994).

\bibitem[DM]{DMorg} P.~Deligne, J.~Morgan. {\em Notes on supersymmetry
  (following Joseph Bernstein)}, in {\em Quantum Fields and Strings: a
  course for mathematicians}, AMS/IAS (1999), 41-97. 


\bibitem[Dr]{DrinfeldQHA} V.~G.~Drinfel'd. {\em Quasi-Hopf
        algebras}, Leningrad J.\ Math.\ 1 (1990) 1419-1457.

\bibitem[Du]{Duflo} M.~Duflo. {\em Caract\`eres des groupes et des
        alg\`ebres de lie r\'esolubles}, Ann.\ Sci.\ ENS 3 (1970),
        23-74.

\bibitem[F]{F} W.~Fulton. {\em Intersection Theory}, Springer (1998).

\bibitem[GM]{GM} S.~I.~Gelfand, Yu.~I.~Manin. {\em Homological
        algebra}, Springer (1999).

\bibitem[GS]{GS} M.~Gerstenhaber, S.~Schack. {\em Algebras,
        bialgebras, quantum groups, and algebraic deformations}, in
        {\em Deformation theory and quantum groups with applications
        to mathematical physics}, Contemp.\ Math.\ 134, AMS (1992),
        51-92.

\bibitem[G]{Guan} D.~Guan. {\em Examples of compact holomorphic
        symplectic manifolds which are not K\"ahlerian, II}, Invent.\
        Math.\ 121 (1995), 135-145.

\bibitem[Ha]{Ha} R.~Hartshorne. {\em Algebraic geometry},
        Springer-Verlag (1977).

\bibitem[HV]{HV} V.~Hinich, A.~Vaintrob. {\em Cyclic operads and
        algebra of chord diagrams}, Selecta Math. 8 (2002), 237-282.

\bibitem[H]{H} N.~Hitchin. {\em Hyper-K\"ahler manifolds}, Ast\'erisque
        206 (1992), 137-166.

\bibitem[HS]{Hitchin-Sawon} N.~Hitchin, J.~Sawon. {\em Curvature and
        characteristic numbers of hyperk\"ahler manifolds}, Duke
        Math.\ J.\ 106 (2001), no. 3, 599-615.

\bibitem[Ka]{Kapranov} M.~Kapranov. {\em Rozansky-Witten
        invariants via Atiyah classes}, Compositio Math.\ 115 (1999),
        71-113.

\bibitem[KS]{KS} M.~Kashiwara, P.~Schapira. {\em Sheaves on manifolds},
        Grundlehren der Mathematischen Wissenschaften 292,
        Springer (1994).

\bibitem[Ks]{Kassel} C.~Kassel. {\em Quantum groups}, Graduate Texts
        in Mathematics 155, Springer (1995).

\bibitem[Ko1]{KontsevichFeynman} M.~Kontsevich. {\em
        Feynman diagrams and low-dimensional topology}, in Proc.\
        First European Congress of Mathematics, Progr.\ Math.\ 120,
        Birkh\"auser, (1994).

\bibitem[Ko2]{KontsDefQuant1} M.~Kontsevich. {\em Deformation
        quantization of Poisson manifolds I}, Lett. Math. Phys.  66
        (2003), 157-216.

\bibitem[Ko3]{KontsDefQuant2} M.~Kontsevich. {\em Operads and motives
        in deformation quantization}, Lett.\ Math.\ Phys.\ 48 (1999),
        35-72.

\bibitem[Ko4]{KontsRW} M.~Kontsevich. {\em Rozansky-Witten invariants
        via formal geometry}, Compositio Math.\ 115 (1999), 115-127.

\bibitem[LMO]{LMO} T.~Q.~T.~Le, J.~Murakami, T.~Ohtsuki. {\em On a
        universal perturbative invariant of $3$-manifolds}, Topology
        37 (1998), 539-574.

\bibitem[LM]{LeMurakamiparallel} T.~Q.~T.~Le,
        J.~Murakami. {\em Parallel version of the universal
        Vassiliev-Kontsevich invariant}, J.\ Pure Appl.\ Algebra 121
        (1997), 271-291.

\bibitem[Ma]{Ma} N.~Markarian. {\em Poincar\'e-Birkhoff-Witt
        isomorphism and Riemann-Roch}, preprint (2000).

\bibitem[MO]{MurTQFT} J.~Murakami, T.~Ohtsuki. {\em Topological quantum
        field theory for the universal quantum invariant}, Comm.\
        Math.\ Phys.\ 188 (1997), 501-520.

\bibitem[NW]{NW} M.~Nieper-Wisskirchen, {\em Chern numbers and
Rozansky-Witten invariants of compact hyper-Kähler manifolds}, World
Scientific (2004).

\bibitem[OG1]{OG1} K.~O'Grady. {\em Desingularized moduli spaces of
        sheaves on a $K3$}, J.\ Reine Angew.\ Math.\ 512 (1999),
        49-117.

\bibitem[OG2]{OG2} K.~O'Grady. {\em A new six dimensional irreducible
        symplectic variety}, J. Algebraic Geom. 12 (2003), 435-505.

\bibitem[Ram]{Ram} A.~Ramadoss. {\em The big Chern classes and the
  Chern character}, arXiv preprint {\tt math.AG/0512104}. 

\bibitem[R]{R} J.~D.~Roberts. {\em Rozansky-Witten theory}, Amer. Math. Soc., Providence, RI, 2003in {\em
Topology and geometry of manifolds (Athens, GA, 2001)}, 1-17,
Proc. Sympos. Pure Math. 71, AMS (2003).



\bibitem[RS]{RS} J.~D.~Roberts, J.~Sawon. {\em Generalisations of
Rozansky-Witten invariants}, in {\em Invariants of knots and
3-manifolds (Kyoto, 2001)}, 263-27, Geom. Topol. Monogr. 4, (2002).

\bibitem[RW]{RW} L.~Rozansky, E.~Witten. {\em
        Hyper-K\"ahler geometry and invariants of $3$-manifolds},
        Selecta Math.\ {\bf 3} (1997) 401--458.

\bibitem[S1]{Sawonthesis} J.~Sawon. {\em Rozansky-Witten
        invariants of hyperk\"ahler manifolds},
        PhD Thesis, University of Cambridge, 1999.

\bibitem[Th]{Th} R.~P.~Thomas. {\em Derived categories for the working
        mathematician}, in {\em Winter School on Mirror Symmetry,
        Vector Bundles and Lagrangian Submanifolds}, 349-361, AMS/IP
        Stud. Adv. Math. 23, AMS (2001). 

\bibitem[Tu]{Tu} V.~G.~Turaev. {\em Quantum invariants of knots and
        3-manifolds}, de Gruyter Studies in Mathematics, 18 (1994).

\bibitem[Va]{Vaintrob} A.~Vaintrob. {\em Vassiliev knot invariants and
        Lie $S$-algebras.}, Math.\ Res.\ Lett.\ 1 (1994), 579-595.

\bibitem[Vo]{Vogel} P.~Vogel. {\em Algebraic structures on modules of
        diagrams}, preprint (1995).

\end{thebibliography}
\end{document}